\input ./style/arxiv-general.cfg
\documentclass[seceqn,aop,MSNbibl,dvips]{arximspdf}
\makeatletter
   \@ifpackageloaded{graphicx}{}{\usepackage{graphicx}}
\makeatother
\usepackage{mathbh}
\usepackage{accents}

\doi{10.1214/14-AOP920} 
\volume{43}
\issue{4}
\pubyear{2015}
\firstpage{1823}
\lastpage{1865}

\makeatletter

\newcommand{\rrVert}{\Vert}
\newcommand{\rrvert}{\vert}
\newcommand{\llVert}{\Vert}
\newcommand{\llvert}{\vert}

\newtheorem{Theorem}{Theorem}[section]
\newtheorem{Lemma}{Lemma}[section]
\newtheorem{Proposition}{Proposition}[section]

\newproclaim{Definition}{Definition}[section]
\newproclaim{Remark}{Remark}[section]

\def\trans{^{\intercal}}

\def\N{\mathbb{N}}
\def\R{\mathbb{R}}
\def\E{\mathbb{E}}
\def\F{\mathbb{F}}
\def\P{\mathbb{P}}

\def\Ac{\mathcal{A}}
\def\Bc{\mathcal{B}}
\def\Ec{\mathcal{E}}
\def\Fc{\mathcal{F}}
\def\Lc{\mathcal{L}}
\def\Pc{\mathcal{P}}
\def\Mc{\mathcal{M}}
\def\Oc{\mathcal{O}}
\def\Vc{\mathcal{V}}

\def\eps{\varepsilon}

\newcommand{\eqref}[1]{(\ref{#1})}
\newcommand{\ess}{\operatorname{ess}}
\newcommand{\tr}{\operatorname{tr}}
\renewcommand{\emptyset}{\varnothing}
\renewcommand{\mathring}[1]{\accentset{\circ}{#1}}
\makeatother

\begin{document}
\begin{frontmatter}

\title{Feynman--Kac representation for
Hamilton--Jacobi--Bellman IPDE}
\runtitle{Feynman--Kac representation for
HJB equation}

\begin{aug}
\author[A]{\fnms{Idris}~\snm{Kharroubi}\corref{}\ead[label=e1]{kharroubi@ceremade.dauphine.fr}}
\and
\author[B]{\fnms{Huy\^en}~\snm{Pham}\ead[label=e2]{pham@math.univ-paris-diderot.fr}}
\runauthor{I. Kharroubi and H. Pham}
\affiliation{Universit\'e Paris Dauphine, and Universit\'e Paris 7
Diderot and CREST-ENSAE}
\address[A]{CEREMADE, CNRS, UMR 7534\\
Universit\'e Paris Dauphine\\
Place Mar\'echal de Lattre de Tassigny\\
75016 Paris\\
France\\
\printead{e1}} 
\address[B]{Laboratoire de Probabilit\'es et\\
\quad Mod\`eles Al\'eatoires, CNRS, UMR 7599\\
Universit\'e Paris 7 Diderot\\
Batiment Sophie Germain, Case 7012\\
Avenue de France\\
75205 Paris Cedex 13\\
France\\
and\\
CREST-ENSAE\\
\printead{e2}}
\end{aug}

\received{\smonth{12} \syear{2012}}
\revised{\smonth{11} \syear{2013}}

%
\begin{abstract}
We aim to provide a Feynman--Kac type representation for
Hamilton--Jacobi--Bellman equation, in terms of forward backward stochastic
differential equation (FBSDE) with a simulatable forward process.
For this purpose, we introduce a class of BSDE where the jumps
component of the solution is subject to a partial nonpositive constraint.
Existence and approximation of a unique minimal solution is proved by a
penalization method under mild assumptions. We then show how minimal
solution to this BSDE class
provides a new probabilistic representation for nonlinear
integro-partial differential equations (IPDEs) of
Hamilton--Jacobi--Bellman (HJB) type, when considering a regime switching forward SDE in a
Markovian framework, and importantly we do not make any ellipticity condition.
Moreover, we state a dual formula of this BSDE minimal solution
involving equivalent change of probability measures.
This gives in particular an original representation for value functions
of stochastic control problems including controlled diffusion coefficient.
\end{abstract}

%

\end{frontmatter}

\section{Introduction}\label{sec1}

The classical Feynman--Kac theorem states that the solution to the
linear parabolic partial differential equation (PDE) of second order:
\begin{eqnarray*}
\frac{\partial v}{\partial t} + b(x).D_x v + \frac{1}{2}\tr\bigl(\sigma
\sigma\trans (x)D_x^2 v\bigr) + f(x) &=& 0,\qquad(t,x)
\in[0,T)\times\R^d,
\\
v(T,x) &=& g(x),\qquad x \in\R^d,
\end{eqnarray*}
may be probabilistically represented under some general conditions as
(see, e.g., \cite{fri75}):
%
%
\begin{eqnarray}
\label{feynlin} v(t,x) &=& \E \biggl[ \int_t^T f
\bigl(X_s^{t,x}\bigr)\,ds + g\bigl(X_T^{t,x}
\bigr) \biggr],
\end{eqnarray}
where $X^{t,x}$ is the solution to the stochastic differential equation
(SDE) driven by a $d$-dimensional Brownian motion $W$ on
a filtered probability space $(\Omega,\Fc,\allowbreak(\Fc_t)_t,\P)$:
\begin{eqnarray*}
dX_s &=& b(X_s)\,ds + \sigma(X_s)
\,dW_s,
\end{eqnarray*}
starting from $x\in\R^d$ at $t\in[0,T]$. By considering the
process $Y_t=v(t,X_t)$, and from It\^o's formula (when $v$ is
smooth) or
in general from martingale representation theorem w.r.t. the Brownian
motion $W$, the Feynman--Kac formula (\ref{feynlin}) is formulated
equivalently in terms of (linear) backward stochastic equation
\begin{eqnarray*}
Y_t &=& g(X_T) + \int_t^T
f(X_s)\,ds - \int_t^T
Z_s \,dW_s,\qquad t \leq T,
\end{eqnarray*}
with $Z$ an adapted process, which is identified to: $Z_t=\sigma
\trans(X_t)D_x v(t,X_t)$ when $v$ is smooth.

Let us now consider the Hamilton--Jacobi--Bellman (HJB) equation in the form:
%
%
\begin{eqnarray}
\label{HJB} &&\frac{\partial v}{\partial t} + \sup_{a \in A} \biggl[
b(x,a).D_x v + \frac{1}{2}\tr\bigl(\sigma \sigma
\trans(x,a)D_x^2 v\bigr) + f(x,a) \biggr]
\nonumber
\\
&&\qquad= 0, \qquad\mbox{on } [0,T)\times\R^d,
\\
&& v(T,x) = g(x),\qquad x \in\R^d,
\nonumber
\end{eqnarray}
where $A$ is a subset of $\R^q$. It is well known (see, e.g., \cite{pha09})
that such nonlinear PDE is the dynamic programming equation
associated to the stochastic control problem with value function
defined by
%
%
\begin{eqnarray}
\label{defsto} v(t,x) &:=& \sup_{\alpha} \E \biggl[ \int
_t^T f\bigl(X_s^{t,x,\alpha
},
\alpha _s\bigr)\,ds + g\bigl(X_T^{t,x,\alpha}\bigr)
\biggr],
\end{eqnarray}
where $X^{t,x,\alpha}$ is the solution to the controlled diffusion
\begin{eqnarray*}
dX_s^\alpha&=& b\bigl(X_s^\alpha,
\alpha_s\bigr)\,ds + \sigma\bigl(X_s^\alpha,
\alpha_s\bigr)\,dW_s,
\end{eqnarray*}
starting from $x$ at $t$, and given a predictable control process
$\alpha$ valued in~$A$.

Our main goal is to provide a probabilistic representation
for the nonlinear HJB equation using backward stochastic differential equations
 (BSDEs), namely the so-called nonlinear Feynman--Kac formula,
which involves a simulatable forward process. One can then hope to use
such representation for deriving a probabilistic numerical scheme for
the solution to HJB equation, hence the stochastic control problem.
Such issues have attracted a lot of interest and generated an important
literature over the recent years. Actually, there is a crucial
distinction between the case where the diffusion coefficient is
controlled or not.

Consider first the case where
$\sigma(x)$ does not depend on $a\in A$, and assume that $\sigma
\sigma\trans(x)$ is of full rank. Denoting by $\theta(x,a)=
\sigma\trans(x)(\sigma\sigma\trans(x))^{-1}b(x,a)$ a solution to
$\sigma(x)\theta(x,a)=b(x,a)$, we notice that the HJB equation
reduces into a
semilinear PDE:
%
%
\begin{eqnarray}
\label{semiF} \frac{\partial v}{\partial t} + \frac{1}{2}\tr\bigl(\sigma\sigma\trans
(x)D_x^2 v\bigr) + F(x,\sigma \trans D_x v)
&=& 0,
\end{eqnarray}
where $F(x,z)=\sup_{a \in A} [f(x,a) + \theta(x,a).z]$ is the
$\theta$-Fenchel--Legendre transform of $f$. In this case, we know from
the seminal works by Pardoux and Peng \cite{parpen90,parpen92},
that the (viscosity) solution $v$ to the semilinear PDE (\ref{semiF})
is connected to the BSDE
%
%
\begin{eqnarray}
\label{BSDEF} Y_t &=& g\bigl(X_T^0\bigr) +
\int_t^T F\bigl(X_s^0,Z_s
\bigr)\,ds - \int_t^T Z_s \,
dW_s,\qquad t \leq T,
\end{eqnarray}
through the relation $Y_t= v(t,X_t^0)$, with a forward diffusion process
\begin{eqnarray*}
dX_s^0 &=& \sigma\bigl(X_s^0
\bigr)\,dW_s.
\end{eqnarray*}
This probabilistic representation leads to a probabilistic numerical
scheme for the resolution to (\ref{semiF}) by discretization and
simulation of
the BSDE (\ref{BSDEF}); see \cite{boutou04}. Alternatively, when the
function $F(x,z)$ is of polynomial type on $z$, the semilinear PDE
(\ref{semiF})
can be numerically solved by a forward Monte--Carlo scheme relying on
marked branching diffusion, as recently pointed out in \cite{hen12}.
Moreover, as showed in \cite{elketal97}, the solution to the BSDE
(\ref{BSDEF}) admits a dual representation in terms of equivalent
change of
probability measures as
%
%
\begin{eqnarray}
\label{dualsto} Y_t &=& \ess\sup_{\alpha}
\E^{\P^\alpha} \biggl[ \int_t^T f
\bigl(X_s^0,\alpha_s\bigr)\,ds + g
\bigl(X_T^0\bigr) \Big| \Fc_t \biggr],
\end{eqnarray}
where for a control $\alpha$, $\P^\alpha$ is the equivalent probability
measure to $\P$ under which
\begin{eqnarray*}
dX_s^0 &=& b\bigl(X_s^0,
\alpha_s\bigr)\,ds + \sigma\bigl(X_s^0\bigr)
\,dW_s^\alpha,
\end{eqnarray*}
with $W^\alpha$ a $\P^\alpha$-Brownian motion by Girsanov's theorem. In
other words,
the process $X^0$ has the same dynamics under $\P^\alpha$ than the
controlled process $X^\alpha$ under $\P$, and the representation
(\ref{dualsto}) can be viewed as a weak formulation (see \cite{elk79})
of the stochastic control problem (\ref{defsto}) in the case of
uncontrolled diffusion coefficient.

The general case with controlled diffusion coefficient $\sigma(x,a)$
associated to fully nonlinear PDE is challenging and led to recent
theoretical advances. Consider the motivating example from uncertain
volatility model in finance formulated here in dimension 1 for
simplicity of notation:
\begin{eqnarray*}
dX_s^\alpha&=& \alpha_s\, dW_s,
\end{eqnarray*}
where the control process $\alpha$ is valued in $A =[\underline
{a},\bar a]$ with $0\leq\underline{a}\leq\bar a<\infty$, and define
the value function of the stochastic control problem
\begin{eqnarray*}
\label{defwvol} v(t,x) & := & \sup_{\alpha} \E\bigl[ g
\bigl(X_T^{t,x,\alpha}\bigr)\bigr],\qquad(t,x) \in [0,T]\times\R.
\end{eqnarray*}
The associated HJB equation takes the form
%
%
\begin{eqnarray}
\label{defGEDP} \frac{\partial v}{\partial t} + G\bigl(D_x^2{v}\bigr)
&=& 0,\qquad(t,x) \in [0,T)\times\R, v(T,x) = g(x), x \in\R,
\end{eqnarray}
where $G(M)=\frac{1}{2}\sup_{a\in A}[a^2 M]=\bar a^2
M^+-\underline{a}^2M^-$. The unique (viscosity) solution to (\ref
{defGEDP}) is represented in terms of the so-called $G$-Brownian motion
$B$, and $G$-expectation $\E_G$, concepts introduced in \cite{pen06}:
\begin{eqnarray*}
v(t,x) &=& \E_G \bigl[g(x+B_{T-t}) \bigr].
\end{eqnarray*}
Moreover, $G$-expectation is closely related to second-order BSDE
studied in \cite{sontouzha11}, namely the process $Y_t =v(t,B_t)$
satisfies a 2BSDE, which is formulated under a nondominated family of
singular probability measures given by the law of $X^\alpha$ under $\P
$. This gives a nice theory and representation for nonlinear PDE, but
it requires a nondegeneracy assumption on the diffusion coefficient,
and does not cover general HJB equation (i.e., control both on drift
and diffusion arising for instance in portfolio optimization). On the
other hand, it is not clear how to simulate $G$-Brownian motion.

We provide here an alternative BSDE representation including general
HJB equation, formulated under a single probability measure (thus
avoiding nondominated singular measures), and under which the forward
process can be simulated. The idea, used in \cite{khaetal10} for
quasi-variational inequalities arising in impulse control problems, is
the following. We introduce a Poisson random measure $\mu_{A}(dt,da)$
on $\R_+\times A$ with finite intensity measure
$\lambda_{A}(da)\,dt$ associated to the marked point process
$(\tau_i,\zeta_i)_i$, independent of $W$, and consider the pure jump
process $(I_t)_t$ equal to the mark $\zeta_i$ valued in $A$ between two
jump times $\tau_i$ and $\tau_{i+1}$. We next consider the forward
regime switching diffusion process
\begin{eqnarray*}
dX_s &=& b(X_s,I_s)\,ds +
\sigma(X_s,I_s)\,dW_s,
\end{eqnarray*}
and observe that the (uncontrolled) pair process $(X,I)$ is Markov. Let
us then consider the BSDE with jumps w.r.t. the
Brownian--Poisson filtration $\F=\F^{W,\mu_{A}}$:
%
%
\begin{eqnarray}
\label{BSDElinjump} Y_t &=& g(X_T) + \int
_t^T f(X_s,I_s)\,ds -
\int_t^T Z_s\, dW_s
\nonumber
\\[-8pt]
\\[-8pt]
&&{}- \int_t^T \int_A
U_s(a) \tilde\mu_A(ds,da),
\nonumber
\end{eqnarray}
where $\tilde\mu_{A}$ is the compensated measure of $\mu_{A}$. This
linear BSDE is the Feynman--Kac formula for the linear
integro-partial differential equation (IPDE):
%
%
\begin{eqnarray}
\label{linIPDE} &&\frac{\partial v}{\partial t} + b(x,a).D_x v +
\frac{1}{2}\tr \bigl(\sigma\sigma\trans (x,a)D_x^2
v\bigr)
\\
&&\qquad{}+ \int_A \bigl(v\bigl(t,x,a'
\bigr)-v(t,x,a)\bigr) \lambda_{A}\bigl(da'\bigr) + f(x,a)
= 0,
\nonumber
\\
\eqntext{\displaystyle(t,x,a) \in[0,T)\times\R^d\times A,}
\\
&&v(T,x,a) = g(x), \qquad(x,a) \in\R^d\times A,
\end{eqnarray}
through the relation: $Y_t=v(t,X_t,I_t)$.
Now, in order to pass from the above linear IPDE with the additional
auxiliary variable $a\in A$ to the nonlinear HJB PDE (\ref{HJB}),
we constrain the jump component to the BSDE (\ref{BSDElinjump}) to be
nonpositive, that is,
%
%
\begin{eqnarray}
\label{consUA} U_t(a) & \leq& 0, \qquad\forall(t,a).
\end{eqnarray}
Then, since $U_t(a)$ represents the jump of $Y_t = v(t,X_t,I_t)$
induced by a jump of the random measure $\mu$, that is of $I$, and assuming
that $v$ is continuous,
the constraint (\ref{consUA}) means that $U_t(a) =
v(t,X_{t},a)-v(t,X_{t},I_{t^-})\leq0$ for all $(t,a)$. This
formally implies that
$v(t,x)$ should not depend on $a \in A$.
Once we get the nondependence of $v$ in $a$, equation (\ref{linIPDE}) becomes
a PDE on $[0,T)\times\R^d$ with a parameter $a \in A$. By taking
the supremum over $a \in A$ in (\ref{linIPDE}), we then obtain
the nonlinear HJB equation (\ref{HJB}).

Inspired by the above discussion, we now introduce the following
general class of BSDE with partially nonpositive jumps, which is a
non-Markovian extension of (\ref{BSDElinjump})--(\ref{consUA}):
%
%
\begin{eqnarray}
\label{BSDEgenintro} Y_t & = & \xi+ \int_t^T
F(s,\omega,Y_s,Z_s,U_s)\,ds +
K_T - K_t
\nonumber
\\[-8pt]
\\[-8pt]
& &{} - \int_t^T Z_s
\,dW_s - \int_t^T\int
_E U_s(e) \tilde\mu (ds,de),\qquad 0 \leq t
\leq T,\mbox{ a.s.}
\nonumber
\end{eqnarray}
with
%
%
\begin{eqnarray}
\label{Uconsintro} U_t(e) & \leq& 0 ,\qquad d\P\otimes dt\otimes
\lambda(de)\qquad \mbox{a.e. on } \Omega\times[0,T]\times A.
\end{eqnarray}
Here, $\mu$ is a Poisson random measure on $\R_+\times E$ with
intensity measure $\lambda(de)\,dt$,
$A$ a subset of $E$, $\xi$ an $\Fc_T$ measurable random variable, and $F$
a generator function. The solution to this BSDE is a quadruple
$(Y,Z,U,K)$ where, besides the usual component $(Y,Z,U)$, the fourth
component $K$ is a predictable nondecreasing process, which makes the
$A$-constraint (\ref{Uconsintro}) feasible. We thus look at the minimal
solution $(Y,Z,U,K)$ in the sense that for any other solution $(\bar
Y,\bar Z,\bar U,\bar K)$ to (\ref{BSDEgenintro})--(\ref{Uconsintro}), we
must have $Y \leq\bar Y$.

We use a penalization method for constructing an approximating sequence
$(Y^n,Z^n,U^n,K^n)_n$
of BSDEs with jumps, and prove that it converges to the minimal
solution that we are looking for. The proof relies on comparison
results, uniform estimates and
monotonic convergence theorem for BSDEs with jumps. Notice that
compared to
\cite{khaetal10}, we do not assume that the intensity measure $\lambda$
of $\mu$ is finite on the whole set $E$, but only on the subset $A$ on
which the jump constraint is imposed.
Moreover in \cite{khaetal10}, the process $I$ does not influence
directly the coefficients of the process $X$, which is Markov in
itself. In contrast, in this paper, we need to enlarge the state
variables by
considering the additional state variable $I$, which makes Markov the
forward regime switching jump-diffusion process $(X,I)$. Our main
result is then
to relate the minimal solution to the BSDE with $A$-nonpositive jumps
to a fully nonlinear IPDE of HJB type
%
%
\begin{eqnarray}
&&\hspace*{-40pt}\frac{\partial v}{\partial t} + \sup_{a\in A} \biggl[
b(x,a). D_x v(t,x) + \frac{1}{2}\tr\bigl(\sigma\sigma\trans
(x,a)D_x^2v(t,x)\bigr)
\nonumber
\\
&&\hspace*{6pt} {}+ \int_{E\setminus A} \bigl[ v\bigl(t,x+\beta(x,a,e)
\bigr) - v(t,x)
\nonumber
\\
&&\hspace*{70pt} {}- \beta(x,a,e).D_x v(t,x) \bigr] \lambda(de)\nonumber
\\
&&\hspace*{92pt} {} + f \bigl(x,a,v,\sigma\trans(x,a)D_x v\bigr)
\biggr] = 0,
\nonumber
\\
\eqntext{\displaystyle\mbox{on } [0,T)\times\R^d.}
\end{eqnarray}
This equation clearly extends HJB equation (\ref{HJB}) by incorporating
integral terms, and with a function $f$ depending on $v$, $D_x v$
(actually, we may also allow
$f$ to depend on integral terms). By the Markov property of the forward
regime switching jump-diffusion process, we easily see that the minimal
solution to the BSDE with $A$-nonpositive jumps is a deterministic
function $v$ of $(t,x,a)$. The main task is to derive the key property
that $v$ does not actually depend on $a$, as a consequence of the
$A$-nonpositive constrained jumps. This issue is a novelty with respect
to the
framework of \cite{khaetal10} where there is a positive cost at each
change of the regime $I$, while in the current paper, the cost is identically
degenerate to zero.
The proof relies on sharp arguments from viscosity solutions,
inf-convolution and semiconcave approximation,
as we do not know a priori any continuity results on $v$.

In the case where the generator function $F$ or $f$ does not depend on
$y,z,u$, which corresponds to the stochastic control framework, we
provide a dual representation of the minimal solution to the BSDE
by means of a family of equivalent change of probability measures in
the spirit of (\ref{dualsto}).
This gives in particular an original representation for value functions
of stochastic control problems, and unifies the weak formulation for
both uncontrolled and controlled diffusion coefficient.

We conclude this introduction by pointing out that our results are
stated without any ellipticity assumption on the diffusion coefficient,
and includes the case of control affecting independently drift and
diffusion, in contrast with the theory of second-order BSDE.
Moreover, our probabilistic BSDE representation leads to a new
numerical scheme for HJB equation, based on the simulation of the
forward process $(X,I)$ and empirical regression methods, hence taking
advantage of the high dimensional properties of Monte--Carlo method.
Convergence analysis for the discrete time approximation of the BSDE
with nonpositive jumps is studied in \cite{KLP13a}, while numerous
numerical tests illustrate the efficiency of the method in \cite{KLP13b}.

The rest of the paper is organized as follows. In Section~\ref{sec2},
we give a
detailed formulation of BSDE with partially nonpositive jumps. We develop
the penalization approach for studying the existence and the
approximation of a unique minimal solution to our BSDE class, and give
a dual representation of the minimal BSDE solution in the stochastic
control case. We show in Section~\ref{sec3} how the minimal BSDE
solution is
related by means of viscosity solutions to the nonlinear IPDE of HJB type.
Finally, we conclude in Section~\ref{sec4} by indicating extensions to our
paper, and discussing probabilistic numerical scheme for the resolution
of HJB equations.

\section{BSDE with partially nonpositive jumps}\label{sec2}


\subsection{Formulation and assumptions}\label{sec2.1}

Let $(\Omega,\Fc,\P)$ be a complete probability space on which are
defined a $d$-dimensional
Brownian motion $W=(W_t)_{t\geq0}$, and an independent integer
valued Poisson random measure $\mu$ on $\R_+\times E$, where $E$ is a
Borelian subset
of $\R^q$, endowed with its Borel $\sigma$-field $\Bc(E)$. We assume
that the
random measure $\mu$ has the
intensity measure $\lambda(de)\,dt$ for some $\sigma$-finite measure
$\lambda$
on $(E, \Bc(E))$ satisfying
\begin{eqnarray*}
\int_E \bigl(1\wedge\llvert e\rrvert ^2
\bigr) \lambda(de) & < & \infty.
\end{eqnarray*}
We set $\tilde\mu(dt,de)=\mu(dt,de)-\lambda(de)\,dt$, the
compensated martingale measure associated to $\mu$, and denote by $\F
=(\Fc_t)_{t\geq0}$ the completion of the natural filtration generated
by $W$ and $\mu$.

We fix a finite time duration $T<\infty$ and we denote by $\Pc$ the
$\sigma$-algebra of $\F$-predictable subsets of $\Omega\times[0,T]$.
Let us
introduce some additional notations.
We denote by:
\begin{itemize}
\item$\mathbf{S^2}$ the set of real-valued c\`adl\`ag adapted processes
$Y=(Y_t)_{0\leq t\leq T}$ such that $ \llVert Y\rrVert_{\mathbf{S^2}}:=
(\E [ \sup_{0\leq t\leq T} \llvert Y_t\rrvert ^2  ] )^{1/2}
< \infty$.

\item$\mathbf{L^p(0,T)}$, $p\geq1$, the set of real-valued
adapted processes $(\phi_t)_{0\leq t\leq T}$ such that $\E [\int_0^T
\llvert \phi_t\rrvert ^p\, dt ]< \infty$.

\item$\mathbf{L^p(W)}$, $p\geq1$, the set of $\R^d$-valued $\Pc
$-measurable processes
$Z=(Z_t)_{0\leq t\leq T}$ such that $\llVert Z\rrVert _{\mathbf{L^p(W)}} :=
(\E [ \int_0^T \llvert Z_t\rrvert ^p \,dt  ] )^{1/ p} <\infty$.

\item$\mathbf{L^p}(\boldsymbol{\tilde\mu})$, $p\geq1$, the set of
$\Pc\otimes\Bc(E)$-measurable maps $U\dvtx\Omega\times[0,T]\times
E\rightarrow\R$ such that $\llVert  U\rrVert _{{\mathbf{L^p}(\boldsymbol
{\tilde\mu})}} :=
(\E [ \int_0^T (\int_E \llvert U_t(e)\rrvert ^2 \lambda(de)
)^{p/2}\,dt
] )^{1/ p} < \infty$.

\item$\mathbf{L^2}(\boldsymbol{\lambda})$ is the set of $\Bc
(E)$-measurable maps\vspace*{-2pt}
$u\dvtx E\rightarrow\R$ such that
$\llvert  u\rrvert _{{\mathbf{L^2}(\boldsymbol{\lambda})}}:=  (\int_E\llvert u(e)\rrvert ^2\lambda
(de)
)^{1/2}< \infty$.\vspace*{2pt}

\item$\mathbf{K^2}$ the closed subset of $\mathbf{S^2}$ consisting of
nondecreasing processes $K= (K_t)_{0\leq t\leq T}$ with $K_0=0$.
\end{itemize}

We are then given three objects:
\begin{enumerate}[3.]
\item[1.] A terminal condition $\xi$, which is an $\Fc_T$-measurable
random variable.
\item[2.] A generator
function $F\dvtx\Omega\times[0,T]\times\R\times\R^d\times
\mathbf{L^2}(\boldsymbol{\lambda})\rightarrow\R$, which is a $\Pc
\otimes\Bc(\R
)\otimes\Bc
(\R^d)\otimes\Bc(\mathbf{L^2}(\boldsymbol{\lambda}))$-measurable map.
\item[3.] A Borelian subset $A$ of $E$ such that $\lambda(A)<\infty$.
\end{enumerate}
We shall impose the following assumption on these objects:

\begin{longlist}[(H0)]
\item[(H0)]
\mbox{}
\begin{enumerate}[(iii)]
\item[(i)] The random variable $\xi$ and the generator function $F$ satisfy
the square integrability condition
\begin{eqnarray*}
\label{xicarint} \E \bigl[\llvert \xi\rrvert ^2 \bigr] + \E \biggl[
\int_{0}^T \bigl\llvert F(t,0,0,0)\bigr\rrvert
^2\,dt \biggr] & < & \infty.
\end{eqnarray*}
%
%
\item[(ii)] The generator function $F$ satisfies the uniform Lipschitz
condition: there exists a constant $C_F$ such that
\begin{eqnarray*}
\bigl\llvert F(t,y,z,u)-F\bigl(t,y',z',u'
\bigr)\bigr\rrvert & \leq& C_F \bigl(\bigl\llvert y-y'
\bigr\rrvert +\bigl\llvert z-z'\bigr\rrvert +\bigl\llvert
u-u'\bigr\rrvert _{{\mathbf{L^2}(\boldsymbol{\lambda})}} \bigr),
\end{eqnarray*}
for all $t \in[0,T]$, $y,y' \in\R$, $z,z'\in\R^d$ and $u,u'\in
\mathbf{L^2}(\boldsymbol{\lambda})$.
%
\item[(iii)] The generator function $F$ satisfies the monotonicity condition
\begin{eqnarray*}
F(t,y,z,u)-F\bigl(t,y,z,u'\bigr) & \leq& \int_E
\gamma \bigl(t,e,y,z,u,u'\bigr) \bigl(u(e)-u'(e)\bigr)
\lambda(de),
\end{eqnarray*}
for all $t \in[0,T]$, $z\in\R^d$, $y \in\R$ and $u,u'\in
\mathbf{L^2}(\boldsymbol{\lambda})$, where $\gamma\dvtx
[0,T]\times\Omega\times
E\times\R
\times\R^d\times\mathbf{L^2}(\boldsymbol{\lambda})\times\mathbf
{L^2}(\boldsymbol{\lambda})\rightarrow\R$
is a $\Pc\otimes\Bc(E)\otimes\Bc(\R)\otimes\Bc(\R^d)\otimes
\Bc(\mathbf{L^2}(\boldsymbol{\lambda}))\otimes\Bc(\mathbf
{L^2}(\boldsymbol{\lambda}))$-measurable map
satisfying:
$C_1(1\wedge\llvert e\rrvert )\leq\gamma(t,e,y,\break z,u,u')\leq C_2(1\wedge
\llvert e\rrvert )$, for all $e \in E$, with two constants $-1<C_1\leq0\leq C_2$.
\end{enumerate}
\end{longlist}


Let us now introduce our class of backward stochastic differential equations
(BSDEs) with partially
nonpositive jumps written in the form
%
%
\begin{eqnarray}
\label{BSDEgen} Y_t & = & \xi+ \int_t^T
F(s,Y_s,Z_s,U_s)\,ds + K_T -
K_t
\nonumber
\\
[-8pt]
\\[-8pt]
& & {}- \int_t^T Z_s
\,dW_s - \int_t^T\int
_E U_s(e) \tilde\mu (ds,de),\qquad 0 \leq t
\leq T,\mbox{ a.s.}
\nonumber
\end{eqnarray}
with
%
%
\begin{eqnarray}
\label{Ucons} U_t(e) & \leq& 0 ,\qquad d\P\otimes dt\otimes
\lambda(de)\qquad \mbox{a.e. on } \Omega\times[0,T]\times A.
\end{eqnarray}

%
\begin{Definition}
A minimal solution to the BSDE with terminal data/\break gener\-ator $(\xi,F)$
and $A$-nonpositive jumps is a quadruple of processes\break $(Y,Z,U,K) \in
\mathbf{S^2}\times\mathbf{L^2(W)}\times\mathbf{L^2}(\boldsymbol
{\tilde\mu})\times\mathbf{K^2}$
satisfying
(\ref{BSDEgen})--(\ref{Ucons})
such that for any other quadruple $(\bar Y,\bar Z,\bar U,\bar K) \in
\mathbf{S^2}\times\mathbf{L^2(W)}\times\mathbf{L^2}(\boldsymbol
{\tilde\mu})\times\mathbf{K^2}$
satisfying (\ref{BSDEgen})--(\ref{Ucons}), we have
\begin{eqnarray*}
Y_t & \leq& \bar Y_t,\qquad0 \leq t \leq T,\mbox{
a.s.}
\end{eqnarray*}
\end{Definition}

%
\begin{Remark}\label{remunicite}
Notice that when it exists, there is a unique minimal solution.
Indeed, by definition, we clearly have uniqueness of the component $Y$.
The uniqueness of $Z$ follows by identifying the Brownian parts and the
finite variation parts, and then the uniqueness of $(U,K)$ is obtained
by identifying the predictable parts and by recalling that the jumps of
$\mu$ are inaccessible.
By misuse of language, we say sometimes that $Y$ [instead of the
quadruple $(Y,Z,U,K)$] is the minimal solution to (\ref
{BSDEgen})--(\ref{Ucons}).
\end{Remark}

In order to ensure that the problem of getting a minimal solution is
well posed, we shall need to assume:

\begin{longlist}
\item[(H1)] There exists a quadruple $(\bar
Y,\bar Z,\bar K,\bar U) \in\mathbf{S^2}\times\mathbf{L^2(W)}\times
\mathbf{L^2}(\boldsymbol{\tilde\mu})\times\mathbf{K^2}$
satisfying (\ref{BSDEgen})--(\ref{Ucons}).
\end{longlist}

We shall see later in Lemma~\ref{lemH1} how such a condition is satisfied
in a Markovian framework.

\subsection{Existence and approximation by penalization}\label{sec2.2}



In this paragraph, we prove the existence of a minimal solution to
(\ref{BSDEgen})--(\ref{Ucons}), based on approximation
via penalization. For each $n \in\N$, we introduce the penalized
BSDE with jumps
%
%
\begin{eqnarray}
\label{BSDEpen} Y_t^n &=& \xi+ \int_t^T
F\bigl(s,Y_s^n,Z_s^n,U_s^n
\bigr)\,ds + K_T^n - K_t^n
\nonumber
\\[-8pt]
\\[-8pt]
& &{} - \int_t^T Z_s^n
\,dW_s - \int_t^T\int
_E U_s^n(e) \tilde\mu(ds,de),
\qquad 0 \leq t \leq T,
\nonumber
\end{eqnarray}
where $K^n$ is the nondecreasing process in $\mathbf{K^2}$ defined by
\begin{eqnarray*}
K_t^n &=& n \int_0^t
\int_A \bigl[U^n_s(e)\bigr]^+
\lambda(de)\,ds,\qquad0 \leq t \leq T.
\end{eqnarray*}
Here, $[u]^+ =\max(u,0)$ denotes the positive part of $u$. Notice
that this penalized BSDE can be rewritten as
\begin{eqnarray*}
Y_t^n &=& \xi+ \int_t^T
F_n\bigl(s,Y_s^n,Z_s^n,U_s^n
\bigr)\,ds - \int_t^T Z_s^n
\,dW_s
\\
&&{}- \int_t^T\int_E
U_s^n(e) \tilde\mu(ds,de),\qquad 0 \leq t \leq T,
\end{eqnarray*}
where the generator $F_n$ is defined by
\begin{eqnarray*}
F_n(t,y,z,u) & = & F(t,y,z,u)+ n\int_A
\bigl[u(e)\bigr]^+\lambda(de),
\end{eqnarray*}
for all $(t,y,z,u)\in[0,T]\times\R\times\R^d\times\mathbf
{L}^2(\lambda
)$. Under {(H0)}(ii)--(iii) and since $\lambda(A)<\infty$, we see
that $F_n$ is Lipschitz continuous w.r.t. $(y,z,u)$ for all $n\in\N$.
Therefore, we obtain from Lemma~2.4 in \cite{tanli94}, that under
(H0), BSDE (\ref{BSDEpen}) admits a unique solution
$(Y^n,Z^n,U^n)\in\mathbf{S^2}\times\mathbf{L^2(W)}\times\mathbf
{L^2}(\boldsymbol{\tilde\mu})$
for any $n\in\N$.



%
\begin{Lemma} \label{leminc} Let Assumption \emph{(H0)} hold.
The sequence $(Y^n)_n$ is nondecreasing, that is, $Y^n_t \leq
Y^{n+1}_t$ for all $ t\in[0, T]$ and all $n \in\N$.
\end{Lemma}

\begin{pf} Fix $n\in\N$, and observe that
\begin{eqnarray*}
F_n(t,e,y,z,u) & \leq& F_{n+1}(t,e,y,z,u),
\end{eqnarray*}
for all $(t,e,y,z,u)\in[0,T]\times E\times\R\times\R^d\times
\mathbf{L^2}(\boldsymbol{\lambda})$.
Under Assumption (H0), we can apply the comparison Theorem~2.5
in \cite{roy06}, which shows that $Y^{n}_t \leq Y^{n+1}_t$,
$0\leq
t\leq T$, a.s.
\end{pf}

The next result shows that the sequence $(Y^n)_n$ is upper-bounded by
any solution to the constrained BSDE.

%
\begin{Lemma} \label{lemcompbor}
Let Assumption \emph{(H0)} hold.
For any quadruple $(\bar Y,\bar Z,\bar U,\bar K) \in\mathbf{
S^2}\times
\mathbf{L^2(W)}\times\mathbf{L^2}(\boldsymbol{\tilde\mu})\times
\mathbf{K^2}$ satisfying
(\ref{BSDEgen})--(\ref{Ucons}), we have
%
%
\begin{eqnarray}
\label{YnleqtildeY} Y_t^n & \leq& \bar Y_t,
\qquad0 \leq t\leq T, n \in\N.
\end{eqnarray}
\end{Lemma}

\begin{pf}
Fix $n\in\N$, and consider a quadruple $(\bar Y,\bar
Z,\bar U,\bar K) \in\mathbf{S^2}\times\mathbf{L^2(W)}\times
\mathbf{L^2}(\boldsymbol{\tilde\mu})\times\mathbf{K^2}$ solution
to (\ref{BSDEgen})--(\ref{Ucons}).
Then, $\bar U$ clearly satisfies $\int_0^t \int_A [\bar
U_s(e)]^+\*\lambda
(de)\,ds = 0$ for all $t \in[0,T]$, and so
$(\bar Y,\bar Z,\bar U,\bar K)$ is a supersolution to the penalized
BSDE (\ref{BSDEpen}), that is,
\begin{eqnarray*}
\bar Y_t &=& \xi+ \int_t^T
F_n(s,\bar Y_s,\bar Z_s,\bar
U_s)\,ds + \bar K_T - \bar K_t
\\
& &{} - \int_t^T \bar Z_s\,
dW_s - \int_t^T\int
_E \bar U_s(e) \tilde \mu(ds,de),\qquad0\leq t
\leq T.
\end{eqnarray*}
By a slight adaptation of the comparison Theorem~2.5 in \cite{roy06}
under {(H0)}, we obtain the required inequality: $Y_t^n \leq
\bar Y_t$, $0\leq t\leq T$.
\end{pf}


We now establish a priori uniform estimates on the sequence
$(Y^n,Z^n,U^n,\allowbreak K^n)_n$.

%
\begin{Lemma} \label{lembor}
Under \emph{(H0)} and \emph{(H1)}, there exists some constant $C$
depending only on $T$ and the monotonicity condition of $F$ in \emph
{(H0)(iii)} such that
%
%
\begin{eqnarray}
\label{bounduni} & &\hspace*{15pt} \bigl\llVert Y^n\bigr\rrVert ^2_{{\mathbf{S^2}}}
+ \bigl\llVert Z^n\bigr\rrVert ^2_{{\mathbf{L^2(W)}}} + \bigl
\llVert U^n\bigr\rrVert ^2_{{\mathbf{L^2}(\boldsymbol{\tilde\mu})}} +\bigl\llVert
K^n\bigr\rrVert ^2_{{\mathbf
{S^2}}}
\nonumber
\\
& &\hspace*{15pt}\qquad\leq C \biggl( \E\llvert \xi\rrvert ^2 + \E \biggl[ \int
_0^T \bigl\llvert F(t,0,0,0)\bigr\rrvert
^2\,dt \biggr] + \E \Bigl[ \sup_{0\leq t\leq T} \llvert \bar
Y_t\rrvert ^2 \Bigr] \biggr),\\
\eqntext{\forall n \in\N.}
\end{eqnarray}
\end{Lemma}

\begin{pf} In what follows, we shall denote by $C>0$ a generic
positive constant depending only on $T$, and the linear growth
condition of $F$ in {(H0)}(ii), which may vary from line to line.
By applying It\^o's formula to $\llvert Y_t^n\rrvert ^2$, and observing that $K^n$ is
continuous and $\Delta Y_t^n = \int_E U_t^n(e)\mu(\{t\},de)$, we have
\begin{eqnarray*}
\E\llvert \xi\rrvert ^2 & = & \E\bigl\llvert Y_t^n
\bigr\rrvert ^2 -2\E\int_t^TY_s^nF
\bigl(s,Y_s^n,Z_s^n,U_s^n
\bigr)\,ds
\\
&&{}- 2\E\int_t^TY_s^n
\,dK^n_s+\E\int_t^T
\bigl\llvert Z^n_s\bigr\rrvert ^2\,ds
\\
& &{} + \E\int_t^T\int_E
\bigl\{\bigl\llvert Y_{s-}^n+U_s^n(e)
\bigr\rrvert ^2-\bigl\llvert Y_{s-}^n\bigr
\rrvert ^2 -2Y_{s-}^nU_s^n(e)
\bigr\}\mu(de,ds)
\\
& = & \E\bigl\llvert Y_t^n\bigr\rrvert ^2 +
\E\int_t^T\bigl\llvert Z^n_s
\bigr\rrvert ^2\,ds + \E\int_t^T
\int_E\bigl\llvert U_s^n(e)\bigr
\rrvert ^2\lambda(de)\,ds
\\
& &{} -2\E\int_t^TY_s^nF
\bigl(s,Y_s^n,Z_s^n,U_s^n
\bigr)\,ds - 2\E\int_t^TY_s^n
\,dK^n_s,\qquad0 \leq t\leq T.
\end{eqnarray*}
From (H0)(iii),
the inequality $Y_t^n \leq\bar{Y}_t$ by Lemma~\ref{lemcompbor} under {(H1)}, and the inequality $2ab \leq
\frac{1}{\alpha}a^2 +\alpha b^2$ for any constant $\alpha>0$,
we have
\begin{eqnarray*}
& & \E\bigl\llvert Y_t^n\bigr\rrvert ^2+\E
\int_t^T\bigl\llvert Z^n_s
\bigr\rrvert ^2\,ds+\E\int_t^T\int
_E\bigl\llvert U_s^n(e)\bigr
\rrvert ^2\lambda(de)\,ds
\\
&&\qquad\leq\E\llvert \xi\rrvert ^2 + C\E\int_t^T
\bigl\llvert Y_s^n\bigr\rrvert \bigl(\bigl\llvert
F(s,0,0,0)\bigr\rrvert +\bigl\llvert Y_s^n\bigr\rrvert +
\bigl\llvert Z_s^n\bigr\rrvert +\bigl\llvert
U_s^n\bigr\rrvert _{\mathbf{L^2}(\boldsymbol{\lambda})} \bigr)\,ds
\\
& &\quad\qquad{}+\frac{1}{\alpha}\E \Bigl[\sup_{s\in[0,T]}\llvert
\bar {Y}_s\rrvert ^2 \Bigr] + \alpha \E\bigl\llvert
K^n_T-K^n_t\bigr\rrvert
^2.
\end{eqnarray*}
Using again the inequality $ab\leq\frac{a^2}{2}+\frac
{b^2}{2}$, and (H0)(i), we get
%
%
\begin{eqnarray}
\label{interYn} & & \E\bigl\llvert Y_t^n\bigr\rrvert
^2+ \frac{1}{2} \E\int_t^T
\bigl\llvert Z^n_s\bigr\rrvert ^2\,ds+
\frac
{1}{2}\E \int_t^T\int
_E\bigl\llvert U_s^n(e)\bigr
\rrvert ^2\lambda(de)\,ds
\nonumber
\\
&&\qquad\leq C \E\int_{t}^T \bigl\llvert
Y_{s}^{n}\bigr\rrvert ^2\,ds + \E\llvert \xi
\rrvert ^2 + \frac
{1}{2}\E\int_0^T
\bigl\llvert F(s,0,0,0)\bigr\rrvert ^2\,ds
\\
&&\quad\qquad{}+\frac{1}{\alpha}\E \Bigl[\sup_{t\in
[0,T]}\llvert
\bar {Y}_t\rrvert ^2 \Bigr] + \alpha\E\bigl\llvert
K^n_T-K^n_t\bigr\rrvert
^2 .
\nonumber
\end{eqnarray}
Now, from the relation (\ref{BSDEpen}), we have\vspace*{1pt}
\begin{eqnarray*}
K^n_T-K^n_t &=&
Y_t^n-\xi-\int_t^T F
\bigl(s,Y_s^n,Z_s^n,U_s^n
\bigr)\,ds
\\
& &{} + \int_t^T Z^n_s
\,dW_s+ \int_t^T\int
_EU_s^n(e)\tilde\mu(ds,de).
\end{eqnarray*}
Thus, there exists some positive constant $C_1$ depending only on the
linear growth condition of $F$ in {(H0)}(ii) such that
%
%
\begin{eqnarray}
\label{inegK}
&&\E\bigl\llvert K^n_T-K^n_t
\bigr\rrvert ^2\nonumber \\
&&\qquad \leq C_1 \biggl( \E\llvert \xi\rrvert
^2 + \E\int_0^T \bigl\llvert
F(s,0,0,0)\bigr\rrvert ^2\,ds + \E\bigl\llvert Y^n_t
\bigr\rrvert ^2
\\
& &\hspace*{65pt}{} + \E\int_t^T \bigl(
\bigl\llvert Y^n_s\bigr\rrvert ^2+\bigl
\llvert Z_s^n\bigr\rrvert ^2+\bigl\llvert
U_s^n\bigr\rrvert _{\mathbf{L^2}(\boldsymbol{\lambda})}^2 \bigr)
\,ds \biggr),\qquad 0 \leq t\leq T.
\nonumber
\end{eqnarray}
Hence, by choosing $\alpha>0$ s.t.
$C_1\alpha\leq\frac{1}{4}$, and plugging into (\ref{interYn}),
we get
\begin{eqnarray*}
& &\frac{3}{4} \E\bigl\llvert Y_t^n\bigr\rrvert
^2+ \frac{1}{4} \E\int_t^T
\bigl\llvert Z^n_s\bigr\rrvert ^2\,ds +
\frac
{1}{4} \E\int_t^T\int
_E\bigl\llvert U_s^n(e)\bigr
\rrvert ^2\lambda(de)\,ds
\\
&&\qquad\leq C \E\int_{t}^T \bigl\llvert
Y_{s}^{n}\bigr\rrvert ^2\,ds +
\frac{5}{4} \E\llvert \xi \rrvert ^2 + \frac
{1}{4} \E\int
_0^T \bigl\llvert F(s,0,0,0)\bigr\rrvert
^2\,ds
\\
&&\quad\qquad{}+ \frac{1}{\alpha} \E \Bigl[\sup_{s\in[0,T]}\llvert
\bar{Y}_s\rrvert ^2 \Bigr],\qquad0 \leq t\leq T.
\end{eqnarray*}
Thus application of Gronwall's lemma to $t\mapsto\E\llvert Y_t^n\rrvert ^2$
yields
%
%
\begin{eqnarray}
\label{interfin} & & \sup_{0\leq t\leq T} \E\bigl\llvert
Y_t^n\bigr\rrvert ^2+ \E\int
_0^T \bigl\llvert Z^n_t
\bigr\rrvert ^2\,dt + \E \int_0^T
\int_E\bigl\llvert U_t^n(e)\bigr
\rrvert ^2\lambda(de)\,dt
\nonumber
\\[-8pt]
\\[-8pt]
&&\qquad\leq C \biggl( \E\llvert \xi\rrvert ^2 + \E\int
_0^T \bigl\llvert F(t,0,0,0)\bigr\rrvert
^2\,dt + \E \Bigl[\sup_{t\in[0,T]}\llvert
\bar{Y}_t\rrvert ^2 \Bigr] \biggr),
\nonumber
\end{eqnarray}
which gives the required uniform estimates (\ref{bounduni}) for
$(Z^n,U^n)_n$ and also $(K^n)_n$ by (\ref{inegK}). Finally, by writing
from (\ref{BSDEpen}) that
\begin{eqnarray*}
\sup_{0\leq t\leq T} \bigl\llvert Y_t^n\bigr
\rrvert & \leq& \llvert \xi\rrvert + \int_0^T
\bigl\llvert F\bigl(t,Y_t^n,Z_t^n,U_t^n
\bigr)\bigr\rrvert \,dt + K_T^n
\\
& &{} + \sup_{0\leq t\leq T} \biggl\llvert \int_0^t
Z^n_s \,dW_s \biggr\rrvert + \sup
_{0\leq t\leq T} \biggl\llvert \int_0^t
\int_E U_s^n(e) \tilde\mu (ds,de)
\biggr\rrvert,
\end{eqnarray*}
we obtain the required uniform estimate (\ref{bounduni}) for $(Y^n)_n$
by the Burkholder--Davis--Gundy inequality, linear growth condition in
{(H0)}(ii), and the uniform estimates for $(Z^n,U^n,K^n)_n$.
\end{pf}

We can now state the main result of this paragraph.

%

%
\begin{Theorem} \label{theomain1}
Under \emph{(H0)} and \emph{(H1)}, there exists a unique minimal solution
$(Y,Z,U,K)\in\mathbf{S^2}\times\mathbf{L^2(W)}\times\mathbf{L^2}
(\boldsymbol{\tilde\mu})\times\mathbf{K^2}$ with $K$ predictable, to
(\ref{BSDEgen})--(\ref{Ucons}). $Y$ is the increasing limit of $(Y^n)_n$
and also in $\mathbf{L^2(0,T)}$,
$K_t$ is the weak limit of $(K^n_t)_n$ in $\mathbf{L^2}
(\boldsymbol{\Omega},\Fc_t,\P)$
for all $t\in[0,T]$, and for any $p\in[1,2)$,
\begin{eqnarray*}
\bigl\llVert Z^n-Z\bigr\rrVert _{{\mathbf{L^p(W)}}} + \bigl\llVert
U^n-U\bigr\rrVert _{{\mathbf
{L^p}(\boldsymbol{\tilde\mu})}} & \longrightarrow& 0,
\end{eqnarray*}
as $n$ goes to infinity.
\end{Theorem}

\begin{pf} By the Lemmata \ref{leminc} and \ref{lemcompbor},
$(Y^n)_n$ converges increasingly to some adapted process $Y$,
satisfying: $\llVert Y\rrVert _{\mathbf{S^2}}<\infty$ by the uniform estimate
for $(Y^n)_n$ in Lemma~\ref{lembor} and Fatou's lemma. Moreover, by
dominated convergence theorem, the convergence of $(Y^n)_n$ to $Y$ also
holds in $\mathbf{L^2(0,T)}$. Next,
by the uniform estimates for $(Z^n,U^n,K^n)_n$ in Lemma~\ref{lembor},
we can apply the monotonic convergence Theorem~3.1 in \cite{esa08},
which extends to the jump case the monotonic convergence theorem of
Peng \cite{pen99} for BSDE. This provides the existence of $(Z,U)\in
\mathbf{L^2(W)}\times\mathbf{L^2}(\boldsymbol{\tilde\mu})$, and
$K$ predictable,
nondecreasing with $\E[K_T^2]<\infty$, such that the sequence
$(Z^n,U^n,K^n)_n$ converges in the sense of Theorem~\ref{theomain1} to
$(Z,U,K)$ satisfying
\begin{eqnarray*}
Y_t &=& \xi+ \int_t^T
F(s,Y_s,Z_s,U_s)\,ds + K_T -
K_t
\\
& & {}- \int_t^T Z_s
\,dW_s - \int_t^T\int
_E U_s(e) \tilde\mu (ds,de),\qquad0 \leq t\leq
T.
\end{eqnarray*}
Thus, the process $Y$ is the difference of a c\`adl\`ag process and
the nondecreasing process $K$, and by Lemma~2.2 in \cite{pen99}, this
implies that $Y$ and $K$ are also c\`adl\`ag, hence respectively
in $\mathbf{S^2}$ and $\mathbf{K^2}$.
Moreover, from the strong convergence in $\mathbf{L^1}(\boldsymbol
{\tilde\mu})$ of
$(U^n)_n$ to $U$ and since $\lambda(A)<\infty$, we have
\begin{eqnarray*}
\E\int_{0}^T\int_A
\bigl[U^{n}_{s}(e)\bigr]^+\lambda(de)\,ds & \longrightarrow
& \E \int_{0}^T\int_A
\bigl[U_{s}(e)\bigr]^+\lambda(de)\,ds,
\end{eqnarray*}
as $n$ goes to infinity. Since $K_T^n=n \int_{0}^T\int_A
[U^{n}_{s}(e)]^+\lambda(de)\,ds$ is bounded in $\mathbf{L^2}(\boldsymbol{\Omega}
,\Fc
_{\mathbf{T}},\P
)$, this implies
\begin{eqnarray*}
\E\int_{0}^T\int_A
\bigl[U_{s}(e)\bigr]^+\lambda(de)\,ds &= & 0,
\end{eqnarray*}
which means that the $A$-nonpositive constraint (\ref{Ucons}) is
satisfied. Hence, $(Y,Z,\allowbreak K,U)$ is a solution to the constrained BSDE
(\ref{BSDEgen})--(\ref{Ucons}), and by\break Lemma~\ref{lemcompbor}, $Y=\lim Y^n$ is the minimal solution. Finally,
the uniqueness of the solution $(Y,Z,U,K)$ is given by Remark~\ref{remunicite}.
\end{pf}

\subsection{Dual representation}\label{sec2.3}

In this subsection, we consider the case where the generator function
$F(t,\omega)$ does not depend on $y,z,u$. Our main goal is to provide a
dual representation of the minimal solution to the BSDE with
$A$-nonpositive jumps in terms of a family of equivalent probability measures.

Let $\Vc$ be the set of $\Pc\otimes\Bc(E)$-measurable processes valued
in $(0,\infty)$, and consider for any $\nu\in\Vc$, the Dol\'
eans--Dade exponential local martingale
%
%
\begin{eqnarray}
\label{doleans} L_t^\nu&:=& \Ec \biggl( \int
_0^{\cdot} \int_E \bigl(
\nu_s(e)-1\bigr)\tilde\mu (ds,de) \biggr)_t
\nonumber
\\
&=& \exp \biggl( \int_0^t \int
_E \ln\nu_s(e)\mu(ds,de) - \int
_0^t \int_E \bigl(
\nu_s(e) - 1\bigr) \lambda(de)\,ds \biggr),
\\
\eqntext{\displaystyle0 \leq t\leq T. }
\end{eqnarray}
When $L^\nu$ is a true martingale, that is, $\E[L_T^\nu]=1$, it
defines a probability measure $\P^\nu$ equivalent to $\P$ on
$(\Omega
,\Fc_T)$ with Radon--Nikodym density
%
%
\begin{eqnarray}
\label{defPnu} \frac{d\P^\nu}{d\P} \bigg|_{\Fc_t} &=& L_t^\nu,
\qquad0 \leq t\leq T,
\end{eqnarray}
and we denote by $\E^\nu$ the expectation operator under $\P^\nu$.
Notice that $W$ remains a Brownian motion under $\P^\nu$, and the
effect of the probability measure $\P^\nu$, by Girsanov's theorem, is
to change the compensator $\lambda(de)\,dt$ of $\mu$ under $\P$ to
$\nu
_t(e)\lambda(de)\,dt$ under
$\P^\nu$. We denote by $\tilde\mu^\nu(dt,de)=\mu(dt,de)-\nu
_t(e)\lambda(de)\,dt$ the compensated martingale measure of $\mu$ under
$\P^\nu$.
We then introduce the subset $\Vc_A$ of $\Vc$ by
\begin{eqnarray*}
\Vc_A &=& \bigl\{ \nu\in\Vc, \mbox{ valued in } [1,\infty) \mbox{
and essentially bounded}\dvtx
\\
& & \hspace*{36.5pt} {}\nu_t(e) = 1, e \in E\setminus A, d\P \otimes
dt\otimes\lambda(de)\mbox{ a.e.} \bigr\},
\end{eqnarray*}
and the subset $\Vc_A^n$ as the elements of $\nu\in\Vc_A$
essentially bounded by $n+1$, for $n\in\N$.

%
\begin{Lemma} \label{lemVCA}
For any $\nu\in\Vc_A$, $L^\nu$ is a uniformly integrable
martingale, and $L_T^\nu$ is square integrable.
\end{Lemma}

\begin{pf}
Several sufficient criteria for $L^\nu$ to be a uniformly integrable
martingale are known. We refer, for example, to the recent paper \cite
{proshi08}, which shows that if
\begin{eqnarray*}
S_T^\nu&:=& \exp \biggl( \int_0^T
\int_E \bigl\llvert \nu_t(e)-1\bigr\rrvert
^2 \lambda (de)\,dt \biggr)
\end{eqnarray*}
is integrable, then $L^\nu$ is uniformly integrable. By definition of
$\Vc_A$, we see that for $\nu\in\Vc_A$,
\begin{eqnarray*}
S_T^\nu&=& \exp \biggl( \int_0^T
\int_A \bigl\llvert \nu_t(e)-1\bigr\rrvert
^2 \lambda(de)\,dt \biggr),
\end{eqnarray*}
which is essentially bounded 
since $\nu$ is essentially bounded and $\lambda(A) <\infty$.
Moreover, from the explicit form (\ref{doleans}) of $L^\nu$, we have
$\llvert L_T^\nu\rrvert ^2=L_T^{\nu^2} S_T^\nu$, and so $\E\llvert L_T^\nu\rrvert ^2\leq\llVert
S_T^\nu\rrVert _\infty$.
\end{pf}

We can then associate to each $\nu\in\Vc_A$ the probability
measure $\P^\nu$ through (\ref{defPnu}).
We first provide a dual representation of the penalized BSDEs in terms
of such $\P^\nu$. To this end, we need the following lemma.

%
\begin{Lemma}\label{lem-vrai-mg}
Let $\phi\in\mathbf{L^2(W)}$ and $\psi\in\mathbf{L^2}(\boldsymbol
{\tilde\mu})$. Then
for every $\nu\in\Vc_A$, the processes $\int_0^{\cdot} \phi_t \,dW_t$ and
$\int_0^{\cdot}\int_E \psi_t(e)\tilde\mu^\nu(dt,de)$ are $\P^\nu$-martingales.
\end{Lemma}

\begin{pf}
Fix $\phi\in\mathbf{L^2(W)}$ and $\nu\in\Vc_A$ and denote
by $M^\phi$ the process $\int_0^{\cdot} \phi_t\,dW_t$. Since $W$ remains a
$\P
^\nu$-Brownian motion, we know that $M^\phi$ is a $\P^\nu$-local
martingale. From the Burkholder--Davis--Gundy and Cauchy--Schwarz inequalites,
we have
\begin{eqnarray*}
\E^\nu \Bigl[\sup_{t\in[0,T]}\bigl\llvert
M^\phi_t\bigr\rrvert \Bigr] & \leq& C\E^\nu
\bigl[\sqrt {\bigl\langle M^\phi\bigr\rangle_{T}} \bigr] = C
\E \biggl[L^\nu_T \sqrt{\int_0^T
\llvert \phi_t\rrvert ^2\,dt} \biggr]
\\
& \leq& C \sqrt{\E \bigl[\bigl\llvert L^\nu_T
\bigr\rrvert ^2 \bigr]}\sqrt{\E \biggl[\int
_0^T\llvert \phi _t\rrvert
^2\,dt \biggr]} < \infty,
\end{eqnarray*}
since $L_T^\nu$ is square integrable by Lemma~\ref{lemVCA}, and $\phi
\in
\mathbf{L^2(W)}$. This implies that $M^\phi$ is $\P^\nu$-uniformly
integrable, and hence a true $\P^\nu$-martingale. The proof for $\int_0^{\cdot}\int_E \phi_t(e)\tilde\mu^\nu(dt,de)$ follows exactly the same
lines and is therefore omitted.
\end{pf}

%
\begin{Proposition} \label{propdualYn}
For all $n\in\N$, the solution to the penalized BSDE (\ref{BSDEpen}) is
explicitly represented as
%
%
\begin{eqnarray}
\label{reldualYn} Y^n_t & =& \operatorname{\ess\sup}\limits
_{\nu\in\Vc_A^{n}}\E
^\nu \biggl[ \xi+\int_t^TF(s)\,ds \Big|
\Fc_t \biggr],\qquad0 \leq t \leq T.
\end{eqnarray}
\end{Proposition}

\begin{pf} Fix $n\in\N$. For any $\nu\in\Vc_A^n$, and by
introducing the compensated martingale measure
$\tilde\mu^\nu(dt,de)=\tilde\mu(dt,de)-(\nu_t(e)-1)\lambda(de)\,dt$
under $\P^\nu$, we see that the solution $(Y^n,Z^n,U^n)$ to the BSDE
(\ref{BSDEpen}) satisfies
%
%
\begin{eqnarray}
\label{Ynnu} Y_t^n &=& \xi+ \int_t^T
\biggl[ F(s) + \int_A \bigl(n\bigl[U_s^n(e)
\bigr]^+ - \bigl(\nu _s(e)-1\bigr) U_s^n(e)
\bigr) \lambda(de) \biggr]\,ds
\nonumber
\\
& & {}- \int_t^T \int_{E\setminus A}
\bigl(\nu_s(e)-1\bigr) U_s^n(e) \lambda
(de)\,ds
\\
&&{}- \int_t^T Z_s^n
\, dW_s - \int_t^T \int
_E U_s^n(e) \tilde\mu
^\nu(ds,de).
\nonumber
\end{eqnarray}
By the definition of $\Vc_A$, we have
\begin{eqnarray*}
\int_t^T \int_{E\setminus A} \bigl(
\nu_s(e)-1\bigr) U_s^n(e) \lambda(de)\,ds &=
& 0,\qquad0\leq t\leq T,\mbox{ a.s.}
\end{eqnarray*}
By taking expectation in (\ref{Ynnu}) under $\P^\nu$ ($\sim\P
$), we
then get from Lemma~\ref{lem-vrai-mg}
%
%
\begin{eqnarray}
\label{Yninter}
\hspace*{18pt}
 Y_t^n &=& \E^\nu \biggl[ \xi+
\int_t^T \biggl( F(s)
\nonumber
\\[-8pt]
\\[-8pt]
\hspace*{18pt}&&\hspace*{60pt} {} + \int_A \bigl(n\bigl[U_s^n(e)
\bigr]^+ - \bigl(\nu_s(e)-1\bigr) U_s^n(e)
\bigr) \lambda(de) \biggr)\,ds \Big| \Fc_t \biggr].
\nonumber
\end{eqnarray}
Now, observe that for any $\nu\in\Vc_A^n$, hence valued in
$[1,n+1]$, we have
\begin{eqnarray*}
n\bigl[U_t^n(e)\bigr]^+ - \bigl(\nu_t(e)-1
\bigr) U_t^n(e) & \geq& 0,\qquad d\P\otimes dt \otimes
\lambda(de) \mbox{ a.e.}
\end{eqnarray*}
which yields by (\ref{Yninter})
%
%
\begin{eqnarray}
\label{Yninter2} Y_t^n & \geq& \operatorname{\ess\sup}
\limits_{\nu\in\Vc
_A^{n}}
\E ^\nu \biggl[ \xi+\int_t^TF(s)\,ds
\Big|\Fc_t \biggr].
\end{eqnarray}
On the other hand, let us consider the process $\nu^*\in\Vc_A^n$
defined by
\begin{eqnarray*}
\nu_t^*(e) &=& \mathbh{1}_{e \in E\setminus A} + \bigl(
\mathbh{1}_{U_t(e)
\leq0} + (n+1) \mathbh{1}_{U_t(e) > 0} \bigr)
\mathbh{1}_{e
\in A},\qquad0 \leq t \leq T, e \in E.
\end{eqnarray*}
By construction, we clearly have
\begin{eqnarray*}
n\bigl[U_t^n(e)\bigr]^+ - \bigl(\nu_t^*(e)-1
\bigr) U_t^n(e) & = & 0,\qquad\forall0 \leq t \leq T, e
\in A,
\end{eqnarray*}
and thus for this choice of $\nu=\nu^*$ in (\ref{Yninter}),
\begin{eqnarray*}
Y_t^n &=& \E^{\nu^*} \biggl[ \xi+\int
_t^TF(s)\,ds \Big|\Fc_t \biggr].
\end{eqnarray*}
Together with (\ref{Yninter2}), this proves the required representation
of $Y^n$.
\end{pf}

%
\begin{Remark}
Arguments in the proof of Proposition~\ref{propdualYn} show that
relation (\ref{reldualYn}) holds for general generator function $F$
depending on
$(y,z,u)$, that is,
\begin{eqnarray*}
Y^n_t & =& \operatorname{\ess\sup}\limits
_{\nu\in\Vc_A^{n}}\E
^\nu \biggl[ \xi+\int_t^TF
\bigl(s,Y_s^n,Z_s^n,U_s^n
\bigr)\,ds \Big|\Fc_t \biggr] ,
\end{eqnarray*}
which is in this case an implicit relation for $Y^n$. Moreover, the
essential supremum in this dual representation is attained for some
$\nu^*$, which takes extreme values $1$ or $n+1$ depending on the sign
of $U^n$, that is, of bang-bang form.
\end{Remark}

Let us then focus on the limiting behavior of the above dual
representation for $Y^n$ when $n$ goes to infinity.

%
\begin{Theorem} \label{theodual}
Under \emph{(H1)}, the minimal solution to (\ref{BSDEgen})--(\ref
{Ucons}) is explicitly represented as
%
%
\begin{eqnarray}
\label{dualY} Y_t & = & \operatorname{\ess\sup}\limits
_{\nu\in\Vc_A}
\E^\nu \biggl[ \xi+\int_t^TF(s)\,ds
\Big|\Fc _t \biggr],\qquad0 \leq t\leq T.
\end{eqnarray}
\end{Theorem}

\begin{pf} Let $(Y,Z,U,K)$ be the minimal solution to (\ref
{BSDEgen})--(\ref{Ucons}).
Let us denote by $\tilde Y$ the process defined in the right-hand side
of (\ref
{dualY}). Since $\Vc_A^n\subset\Vc_A$, it is clear from
the representation (\ref{reldualYn}) that $Y_t^n \leq\tilde Y_t$,
for all $n$. Recalling from Theorem~\ref{theomain1} that $Y$ is the
pointwise limit of $Y^n$, we deduce that $Y_t = \lim_{n\rightarrow
\infty} Y_t^n \leq\tilde Y_t$, \mbox{$0 \leq t\leq T$}.

Conversely, for any $\nu\in\Vc_A$, let us consider the
compensated martingale measure
$\tilde\mu^\nu(dt,de) = \tilde\mu(dt,de)-(\nu_t(e)-1)\lambda
(de)\,dt$
under $\P^\nu$, and observe that $(Y,Z,U,K)$ satisfies
%
%
\begin{eqnarray}
\label{Ynu} Y_t &=& \xi+ \int_t^T
\biggl[ F(s) - \int_A \bigl(\nu_s(e)-1\bigr)
U_s(e) \lambda (de) \biggr]\,ds + K_T - K_t
\nonumber
\\
& &{} - \int_t^T \int_{E\setminus A}
\bigl(\nu_s(e)-1\bigr) U_s(e) \lambda (de)\,ds
\\
&&{} - \int_t^T Z_s \,
dW_s - \int_t^T \int
_E U_s(e) \tilde\mu^\nu(ds,de).
\nonumber
\end{eqnarray}
By the definition of $\nu\in\Vc_A$, we have $\int_t^T \int_{E\setminus A} (\nu_s(e)-1) U_s(e) \lambda(de)\,ds = 0$. Thus,\vspace*{1pt} by
taking expectation in (\ref{Ynu}) under $\P^\nu$ from Lemma~\ref
{lem-vrai-mg}, and recalling that $K$ is nondecreasing, we have
\begin{eqnarray*}
Y_t & \geq& \E^\nu \biggl[ \xi+ \int
_t^T \biggl( F(s) - \int_A
\bigl(\nu _s(e)-1\bigr) U_s(e) \lambda(de) \biggr)\,ds
\Big| \Fc_t \biggr]
\\
& \geq& \E^\nu \biggl[ \xi+ \int_t^T
F(s)\,ds \Big| \Fc_t \biggr],
\end{eqnarray*}
since $\nu$ is valued in $[1,\infty)$, and $U$ satisfies the
nonpositive constraint (\ref{Ucons}). Since $\nu$ is arbitrary in
$\Vc
_A$, this proves
the inequality $Y_t \geq\tilde Y_t$, and finally the required
relation $Y = \tilde Y$.
\end{pf}

\section{Nonlinear IPDE and Feynman--Kac formula}\label{sec3}

In this section, we shall show how minimal solutions to our BSDE class
with partially nonpositive jumps provides actually a new probabilistic
representation (or the Feynman--Kac formula) to fully nonlinear
integro-partial differential equation (IPDE) of
Hamilton--Jacobi--Bellman (HJB) type,
when dealing with a suitable\break Markovian framework.

\subsection{The Markovian framework}\label{sec3.1}

We are given a compact set $A$ of $\R^q$, and a Borelian subset $L
\subset\R^l\setminus\{0\}$, equipped with respective Borel
$\sigma
$-fields $\Bc(A)$ and $\Bc(L)$. We assume that:

\begin{longlist}[(H$A$)]
\item[(H$A$)] The interior set $\mathring{A}$ of $A$
is connex, and $A = \operatorname{Adh}(\mathring{A})$, the closure of
its interior.
\end{longlist}

We consider the case where $E = L\cup A$ and we may assume w.l.o.g.
that $L\cap A=\emptyset$ by identifying $A$ and $L$,
respectively, with the sets $A\times\{0\}$ and $\{0\}\times L$ in $\R
^q\times\R^l$. We consider two independent Poisson random measures
$\vartheta$ and $\pi$ defined respectively on $\R_+\times L$ and $\R
_+\times A$. We suppose that $\vartheta$ and $\pi$ have respective
intensity measures
$\lambda_\vartheta(d\ell)\,dt$ and $\lambda_\pi(da)\,dt$ where
$\lambda
_\vartheta$ and $\lambda_\pi$ are two $\sigma$-finite measures with
respective supports $L$ and $A$, and satisfying
\begin{eqnarray*}
\int_L\bigl(1\wedge\llvert \ell\rrvert ^2
\bigr)\lambda_\vartheta(d\ell) < \infty \quad \mbox{and}\quad\int
_A\lambda_\pi(da) < \infty,
\end{eqnarray*}
and we denote by $\tilde\vartheta(dt,d\ell) = \vartheta
(dt,d\ell
)-\lambda_\vartheta(d\ell)\,dt$ and $\tilde\pi(dt,da) = \pi
(dt,da)-\lambda_\pi(da)\,dt$ the compensated martingale measures of
$\vartheta$ and $\pi$, respectively. We also assume that:

\begin{longlist}[(H$\lambda_\pi$)]
\item[(H$\lambda_\pi$)]
\mbox{}
\begin{enumerate}[(ii)]
\item[(i)] The measure $\lambda_\pi$ supports the whole set
$\mathring A$:
for any $a \in\mathring A$ and any open neighborhood $\Oc$ of $a$
in $\R^q$
we have $\lambda_\pi(\Oc\cap\mathring A) > 0$.
\item[(ii)] The boundary of $A$: $\partial A = A\setminus\mathring
A$, is
negligible w.r.t. $\lambda_\pi$, that is, $\lambda_\pi(\partial A)
= 0$.
\end{enumerate}
\end{longlist}

In this context, by taking a random measure $\mu$ on
$\R_+\times E$ in the form, $\mu= \vartheta+\pi$, we notice that
it remains a Poisson random measure with intensity measure
$\lambda(de)\,dt$ given by
\begin{eqnarray*}
\int_E\varphi(e)\lambda(de) & = & \int
_L\varphi(\ell)\lambda _\vartheta (d\ell)+\int
_A\varphi(a)\lambda_\pi(da),
\end{eqnarray*}
for any measurable function $\varphi$ from $E$ to $\R$, and
we have the following identifications:
%
%
\begin{eqnarray}
\label{identif L2mu} \mathbf{L^2}(\boldsymbol{\tilde\mu}) =
\mathbf{L^2}(\boldsymbol {\tilde\vartheta})\times
\mathbf{L^2}(\boldsymbol{\tilde\pi}), \qquad 
\mathbf{L^2}(\boldsymbol{\lambda}) =
\mathbf{L^2}(\boldsymbol{\lambda_ \vartheta})\times \mathbf{L^2}(\boldsymbol{\lambda_\pi}),
\end{eqnarray}
where:
\begin{itemize}
\item$ \mathbf{L^2}(\boldsymbol{\tilde\vartheta})$ is the set of
$\Pc
\otimes
\Bc(L)$-measurable maps $U\dvtx\Omega\times[0,T]\times L\rightarrow
\R$
such that
\begin{eqnarray*}
{\llVert U\rrVert }_{\mathbf{L^2}(\boldsymbol{\tilde\vartheta})} & := & \biggl(\E \biggl[\int
_0^T\int_L\bigl\llvert
U_t(\ell)\bigr\rrvert ^2\lambda_\vartheta(d\ell)
\,dt \biggr] \biggr)^{1/2}<\infty.
\end{eqnarray*}
\item$ \mathbf{L^2}(\boldsymbol{\tilde\pi})$ is the set of $\Pc
\otimes\Bc
(A)$-measurable maps $R\dvtx\Omega\times[0,T]\times A\rightarrow\R$
such that
\begin{eqnarray*}
{\llVert R\rrVert }_{\mathbf{L^2}(\boldsymbol{\tilde\pi})} & := & \biggl(\E \biggl[\int
_0^T\int_A \bigl\llvert
R_t(a)\bigr\rrvert ^2\lambda_\pi(da)\,dt
\biggr] \biggr)^{1/2}<\infty.
\end{eqnarray*}
\item$ \mathbf{L^2}(\boldsymbol{\lambda_\vartheta})$ is the set of
$\Bc
(L)$-measurable maps $u\dvtx L\rightarrow\R$ such that
\begin{eqnarray*}
{\llvert u\rrvert }_{\mathbf{L^2}(\boldsymbol{\lambda_\vartheta})} & := & \biggl(\int_L
\bigl\llvert u(\ell )\bigr\rrvert ^2\lambda_\vartheta(d\ell)
\biggr)^{1/2}<\infty.
\end{eqnarray*}
\item$\mathbf{L^2}(\boldsymbol{\lambda_\pi})$ is the set of $\Bc
(A)$-measurable
maps $r\dvtx A\rightarrow\R$ such that
\begin{eqnarray*}
{\llvert r\rrvert }_{\mathbf{L^2}(\boldsymbol{\lambda_\pi})} & := & \biggl(\int_A
\bigl\llvert r(a)\bigr\rrvert ^2\lambda_\pi (da)
\biggr)^{1/2}<\infty.
\end{eqnarray*}
\end{itemize}

Given some measurable functions $b \dvtx \R^d\times\R^q \rightarrow
\R^d$, $\sigma\dvtx \R^d\times\R^q \rightarrow\R^{d\times d}$
and $\beta\dvtx \R^d\times\R^q\times L \rightarrow\R^d$, we
introduce the forward Markov regime-switching jump-diffusion process
$(X,I)$ governed by
%
%
\begin{eqnarray}
\label{dynX}\hspace*{8pt} dX_s &=& b(X_s,I_s)\,ds +
\sigma(X_s,I_s)\,dW_s + \int
_L \beta (X_{s^-},I_{s^-},\ell) \tilde
\vartheta(ds,d\ell),
\\
\label{dynI}\hspace*{8pt} dI_s &=& \int_A
(a-I_{s^-} ) \pi(ds,da).
\end{eqnarray}
In other words, $I$ is the pure jump process valued in $A$ associated
to the Poisson random measure $\pi$, which changes the coefficients of
jump-diffusion process $X$. We make the usual assumptions on the
forward jump-diffusion coefficients:
\begin{longlist}
\item[(HFC)]
\mbox{}
\begin{enumerate}[(ii)]
\item[(i)] There exists a constant $C$ such that
\begin{eqnarray*}
\bigl\llvert b(x,a)-b\bigl(x',a'\bigr)\bigr\rrvert +
\bigl\llvert \sigma(x,a)-\sigma\bigl(x',a'\bigr)\bigr
\rrvert & \leq& C \bigl(\bigl\llvert x-x'\bigr\rrvert +\bigl\llvert
a-a'\bigr\rrvert \bigr),
\end{eqnarray*}
for all $x,x'\in\R^d$ and $a,a'\in\R^q$.
\item[(ii)] There exists a constant $C$ such that
%
\begin{eqnarray*}
\bigl\llvert \beta(x,a,\ell)\bigr\rrvert & \leq& C\bigl(1+\llvert x\rrvert
\bigr) \bigl(1 \wedge\llvert \ell \rrvert \bigr),
\\
\bigl\llvert \beta(x,a,\ell) - \beta\bigl(x',a',\ell
\bigr)\bigr\rrvert & \leq& C \bigl(\bigl\llvert x-x'\bigr\rrvert +
\bigl\llvert a-a'\bigr\rrvert \bigr) \bigl(1 \wedge\llvert \ell
\rrvert \bigr),
\end{eqnarray*}
for all $x,x'\in\R^d$, $a,a' \in\R^q$ and $\ell\in L$.
\end{enumerate}
\end{longlist}

%
\begin{Remark}
We do not make any ellipticity assumption on $\sigma$. In
particular, some lines and columns of $\sigma$ may be equal to zero,
and so there is no loss of generality by considering that the dimension
of $X$ and $W$ are equal. We can also make the coefficients $b,\sigma$
and $\beta$ depend on time with the following standard procedure: we
introduce the time variable as a state component $\Theta_t = t$,
and consider the forward Markov system:
\begin{eqnarray*}
\hspace*{1pt}dX_s &=& b(X_s,\Theta_s,I_s)
\,ds + \sigma(X_s,\Theta_s,I_s)
\,dW_s + \int_L \beta(X_{s^-},
\Theta_{s^-},I_{s^-},\ell) \tilde\vartheta (ds,d\ell),
\\
\hspace*{1pt}d\Theta_s &=& ds,
\\
\hspace*{1pt}dI_s &=& \int_A (a-I_{s^-} )
\pi(ds,da),
\end{eqnarray*}
which is of the form given above, but with an enlarged state $(X,\Theta
,I)$ (with degenerate noise), and with the resulting assumptions on
$b(x,\theta,a)$, $\sigma(x,\theta,a)$ and $\beta(x,\theta,a,\ell)$.
\end{Remark}

Under these conditions, existence and uniqueness of a solution
$(X^{t,x,a}_s,\allowbreak I_s^{t,a})_{t\leq s\leq T}$ to (\ref{dynX})--(\ref{dynI})
starting from
$(x,a) \in\R^d\times\R^q$ at time $s=t \in[0,T]$, is well
known, and we have the standard estimate: for all $p \geq2$,
there exists some positive constant $C_p$ s.t.
%
%
\begin{eqnarray}
\label{estimX} \E \Bigl[ \sup_{t\leq s\leq T} \bigl\llvert
X_s^{t,x,a} \bigr\rrvert ^p + \bigl\llvert
I_s^{t,a}\bigr\rrvert ^p \Bigr] & \leq&
C_p\bigl(1 + \llvert x\rrvert ^p+\llvert a\rrvert
^p\bigr),
\end{eqnarray}
for all $(t,x,a) \in[0,T]\times\R^d\times\R^q$.

In this Markovian framework, the terminal data and generator of our
class of BSDE are given by
two continuous functions $g\dvtx \R^d\times\R^q\rightarrow\R$ and
$f \dvtx \R^d\times\R^q\times\R\times\R^d\times\mathbf
{L^2}(\boldsymbol{\lambda
_\vartheta})\rightarrow\R$. We make the following assumptions on the
BSDE coefficients:
\begin{longlist}
\item[(HBC1)]
\mbox{}
\begin{enumerate}[(ii)]
\item[(i)] The functions $g$ and $f(\cdot,0,0,0)$ satisfy a
polynomial growth condition:
\begin{eqnarray*}
\sup_{x\in\R^d, a\in\R^q} \frac{ \llvert g(x,a)\rrvert  + \llvert f(x,a,0,0,0)\rrvert }{1 +
\llvert x\rrvert ^m+\llvert a\rrvert ^m} & < & \infty,
\end{eqnarray*}
for some $m \geq0$.
\item[(ii)] There exists some constant $C$ s.t.
%
\begin{eqnarray*}
& & \bigl\llvert f(x,a,y,z,u) - f\bigl(x',a',y',z',u'
\bigr)\bigr\rrvert
\\
&&\qquad \leq C \bigl(\bigl\llvert x-x'\bigr\rrvert +\bigl\llvert
a-a'\bigr\rrvert +\bigl\llvert y-y'\bigr\rrvert +
\bigl\llvert z-z'\bigr\rrvert + {\bigl\llvert u-u'
\bigr\rrvert }_{\mathbf{L^2}(\boldsymbol{\lambda_\vartheta})} \bigr),
\end{eqnarray*}
for all $x,x'\in\R^d$, $y,y' \in\R$, $z,z' \in\R^d$, $a,a'\in\R^q$
and $u,u'\in\mathbf{L^2}(\boldsymbol{\lambda_\vartheta})$.
\end{enumerate}
\end{longlist}

\begin{longlist}
\item[(HBC2)]
The generator function $f$ satisfies the monotonicity condition:
\begin{eqnarray*}
&&f(x,a,y,z,u)-f\bigl(x,a,y,z,u'\bigr)
\\
&&\qquad\leq\int_L \gamma\bigl(x,a,\ell,y,z,u,u'\bigr) \bigl(u(\ell)-u'(\ell)\bigr)
\lambda_\vartheta(d\ell),
\end{eqnarray*}
for all $x \in\R^d$, $a\in\R^q$, $z\in\R^d$, $y \in\R$ and
$u,u'\in
\mathbf{L^2}(\boldsymbol{\lambda_\vartheta})$, where $\gamma\dvtx \R^d\times
E\times\R\times\R^d\times\mathbf{L^2}(\boldsymbol{\lambda_\vartheta})\times\mathbf{L^2}(\boldsymbol{\lambda_\vartheta})\rightarrow\R$ is a $\Bc
(\R
^d)\otimes\Bc(E)\otimes\Bc(\R)\otimes\Bc(\R^d)\otimes\Bc
(\mathbf
{L^2}(\boldsymbol{\lambda_\vartheta}))\otimes\Bc(\mathbf{L^2}(\boldsymbol{\lambda_\vartheta}))$-measurable map
satisfying: $C_1(1\wedge\llvert \ell\rrvert ) \leq\gamma(x,a,\ell,y,\break z,u,u') \leq
C_2(1\wedge\llvert \ell\rrvert )$, for $\ell\in L$, with two constants
$-1<C_1\leq0\leq C_2$.
\end{longlist}

Let us also consider an assumption on the dependence of $f$ w.r.t. the
jump component used in \cite{barbucpar97}, and stronger than
(HBC2).
\begin{longlist}
\item[(HBC2$'$)]
 The generator function $f$ is of the form
\begin{eqnarray*}
f(x,a,y,z,u) & = & h \biggl(x,a,y,z,\int_Lu(\ell)
\delta(x,\ell )\lambda _\vartheta(d\ell) \biggr)
\end{eqnarray*}
for $(x,a,y,z,u)\in\R^d\times\R^q\times\R\times\R^d\times
\mathbf{L^2}(\boldsymbol{\lambda})$, where:
\begin{itemize}
\item$\delta$ is a measurable function on $\R^d\times L$ satisfying:
\begin{eqnarray*}
0 &\leq&\delta(x,\ell) \leq C\bigl(1\wedge\llvert \ell\rrvert \bigr),
\\
\bigl\llvert \delta(x,\ell)-\delta\bigl(x',\ell\bigr)\bigr\rrvert
& \leq& C\bigl\llvert x-x'\bigr\rrvert \bigl(1\wedge\llvert \ell
\rrvert ^2\bigr),\qquad x,x' \in\R^d, \ell
\in L,
\end{eqnarray*}
for some positive constant $C$.
\item
$h$ is a continuous function on $\R^d\times\R^q\times\R\times\R
^d\times
\R$ such that $\rho\mapsto h(x,a,y,z,\rho)$
is nondecreasing for all $(x,a,y,z)\in\R^d\times\R^q\times\R
\times\R
^d$, and satisfying for some positive constant $C$:
\begin{eqnarray*}
\bigl\llvert h(x,a,y,z,\rho)-h\bigl(x,a,y,z,\rho'\bigr)\bigr
\rrvert & \leq& C\bigl\llvert \rho-\rho'\bigr\rrvert,\qquad \rho,
\rho' \in\R,
\end{eqnarray*}
for all $(x,a,y,z) \in\R^d\times\R^q\times\R\times\R^d$.
\end{itemize}
\end{longlist}

Now with the identification (\ref{identif L2mu}), the BSDE problem
(\ref{BSDEgen})--(\ref{Ucons}) takes the following form:
find the minimal solution $(Y,Z,U,R,K) \in\mathbf{S^2}\times\mathbf
{L^2(W)}\times\mathbf{L^2}(\boldsymbol{\tilde\vartheta})\times\mathbf{L^2}(\boldsymbol{\tilde\pi})\times\mathbf{K^2}$ to
%
%
\begin{eqnarray}
\label{BSDEmarkov}
\hspace*{9pt}
Y_t &=& g(X_T,I_T) + \int
_t^T f (X_s,I_s,Y_s,Z_s,U_s
)\,ds + K_T - K_t
\nonumber
\\[-8pt]
\\[-8pt]
\hspace*{9pt}& & {}- \int_t^T Z_s.dW_s
- \int_t^T \int_L
U_s(\ell) \tilde \vartheta(ds,d\ell)-\int_t^T
\int_A R_s(a) \tilde\pi (ds,da),
\nonumber
\end{eqnarray}
with
%
%
\begin{eqnarray}
\label{cons-markov} R_t(a) & \leq& 0 ,\qquad d\P\otimes dt\otimes
\lambda_\pi(da)\mbox{ a.e.}
\end{eqnarray}

The main goal of this paper is to relate the BSDE (\ref{BSDEmarkov})
with $A$-nonpositive jumps (\ref{cons-markov}) to the following
nonlinear IPDE of HJB type:
%
%
\begin{eqnarray}
\label{IPDE} && {-} \frac{\partial w}{\partial t} - \sup_{a\in A} \bigl[
\Lc^a w + f \bigl(\cdot,a,w,\sigma\trans (\cdot,a)D_x w,
\Mc^a w\bigr) \bigr]= 0,
\\
\eqntext{\mbox{on } [0,T)\times\R^d,}
\\
\label{termIPDE} && w(T,x) = \sup_{a\in A} g(x,a),\qquad x \in
\R^d,
\end{eqnarray}
where
\begin{eqnarray*}
\Lc^a w(t,x) &=& b(x,a). D_x w(t,x) +
\frac{1}{2}\tr\bigl(\sigma \sigma \trans(x,a)D_x^2w(t,x)
\bigr)
\\
& & {}+ \int_L \bigl[ w\bigl(t,x+\beta(x,a,\ell)\bigr) -
w(t,x)
\\
&&\hspace*{56pt} {}- \beta (x,a,\ell).D_x w(t,x) \bigr]
\lambda_\vartheta(d\ell),
\\
\Mc^a w(t,x) &=& \bigl(w\bigl(t,x+\beta(x,a,\ell)\bigr)- w(t,x)
\bigr)_{\ell
\in L},
\end{eqnarray*}
for $(t,x,a)\in[0,T]\times\R^d\times\R^q$.

Notice that under (HBC1), (HBC2) and (\ref{estimX})
[which follows from {(HFC)}], and with the identification (\ref
{identif L2mu}),
the generator $F(t,\omega,y,z,u,r) = f(X_t(\omega),\break I_t(\omega
),y,\allowbreak z,u)$ and the terminal condition $\xi= g(X_T,I_T)$ satisfy
clearly Assumption (H0). Let us now show that Assumption {(H1)}
is satisfied. More precisely, we have the following result.

%
\begin{Lemma} \label{lemH1}
Let Assumptions \emph{(HFC)}, \emph{(HBC1)} hold.
Then, for any initial condition $(t,x,a)\in[0,T]\times\R^d\times\R^q$,
there exists a solution
$\{(\bar Y^{t,x,a}_s,\bar Z ^{t,x,a}_s,\break \bar U^{t,x,a}_s,\allowbreak\bar R
^{t,x,a}_s,\bar K ^{t,x,a}_s),t\leq s\leq T\}$ to the BSDE
(\ref{BSDEmarkov})--(\ref{cons-markov}) when $(X,I) =\break \{
(X_s^{t,x,a},I_s^{t,a})$, $t \leq s\leq T\}$, with
$\bar Y_s^{t,x,a} = \bar v(s,X_s^{t,x,a})$ for some deterministic
function $\bar v$ on $[0,T]\times\R^d$ satisfying a polynomial growth
condition:
for some $p \geq2$,
%
%
\begin{eqnarray}
\label{polybarv} \sup_{(t,x)\in[0,T]\times\R^d}\frac{\llvert \bar v(t,x)\rrvert }{1+\llvert x\rrvert ^p} & < & \infty.
\end{eqnarray}
\end{Lemma}

\begin{pf} Under {(HBC1}) and since $A$ is compact, we observe
that there exists some $m \geq0$ such that
%
%
\begin{eqnarray}
\label{growthgf} C_{f,g} := \sup_{x\in\R^d,a\in A}
\frac{\llvert g(x,a)\rrvert  +
\llvert f(x,a,y,z,u)\rrvert }{1 + \llvert x\rrvert ^m + \llvert y\rrvert  + \llvert z\rrvert  + \llvert u\rrvert _{\mathbf{L^2}
(\boldsymbol{\lambda_\vartheta})}} & < & \infty.
\end{eqnarray}
Let us then consider the smooth function $\bar v(t,x) = \bar C
e^{\rho(T-t)}(1+\llvert x\rrvert ^p)$ for some positive constants $\bar C$ and $\rho$
to be determined later, and with $p = \max(2,m)$.
We claim that for $\bar C$ and $\rho$ large enough, the function $\bar
v$ is a classical supersolution to (\ref{IPDE})--(\ref{termIPDE}).
Indeed, observe first that from the growth condition on $g$ in (\ref
{growthgf}), there exists $\bar C > 0$ s.t.
$\hat g(x) := \sup_{a\in A} g(x,a) \leq\bar C(1+\llvert x\rrvert ^p)$ for
all $x \in\R^d$. For such $\bar C$, we then have: $\bar v(T,\cdot)
\geq
\hat g$.
On the other hand, we see after straightforward calculation that there
exists a positive constant $C$ depending only on $\bar C$, $C_{f,g}$, and
the linear growth condition in $x$ on $b$, $\sigma$, $\beta$ by
{(HFC)} (recall that $A$ is compact), such that
\begin{eqnarray*}
- \frac{\partial\bar v}{\partial t} - \sup_{a\in A} \bigl[ \Lc^a
\bar v + f \bigl(\cdot,a,\bar v,\sigma\trans(\cdot,a)D_x \bar v,
\Mc^a \bar v\bigr) \bigr] & \geq& (\rho- C) \bar v
\\
& \geq& 0,
\end{eqnarray*}
by choosing $\rho\geq C$. Let us now define the quintuple $(\bar
Y,\bar Z,\bar U,\bar R,\bar K)$ by
\begin{eqnarray*}
\bar Y_t &=& \bar v(t,X_t) \qquad\mbox{for } t < T,
\bar Y_T = g(X_T,I_T),
\\
\bar Z_t &=& \sigma\trans(X_{t^-},I_{t^-})D_x
\bar v(t,X_{t^-}),\qquad t \leq T,
\\
\bar U_t &=& \Mc^{I_{t^-}} \bar v(t,X_{t^-}), \qquad
\bar R_t = 0, t \leq T ,
\\
\bar K_t &=& \int_0^t \biggl[ -
\frac{\partial\bar v}{\partial
t}(s,X_s) - \Lc^{I_s}\bar
v(s,X_s) - f(X_s,I_s,\bar
Z_s,\bar U_s) \biggr]\,ds,\qquad t < T,
\\
\bar K_T &=& \bar K_{T^-} + \bar v(T,X_T) -
g(X_T,I_T).
\end{eqnarray*}
From the supersolution property of $\bar v$ to (\ref{IPDE})--(\ref
{termIPDE}), the process $\bar K$ is nondecreasing.
Moreover, from the polynomial growth condition on $\bar v$, linear
growth condition on $b$, $\sigma$, growth condition (\ref{growthgf}) on
$f$, $g$
and the estimate (\ref{estimX}), we
see that $(\bar Y,\bar Z,\bar U,\bar R,\bar K)$ lies in
$\mathbf{S^2}\times\mathbf{L^2(W)}\times\mathbf{L^2}(\boldsymbol
{\tilde{\vartheta}})\times\mathbf{L^2}(\boldsymbol{\tilde{\pi
}})\times\mathbf{K^2}$.
Finally, by applying It\^o's formula to $\bar v(t,X_t)$, we conclude
that $(\bar Y,\bar Z,\bar U,\bar R,\bar K)$ is solution a to (\ref
{BSDEmarkov}), and the constraint (\ref{cons-markov}) is trivially satisfied.
\end{pf}

Under (HFC), (HBC1) and (HBC2), we then get
from Theorem~\ref{theomain1} the existence of a unique minimal solution
$\{(Y_s^{t,x,a},Z_s^{t,x,a},U_s^{t,x,a},R_s^{t,x,a},K_s^{t,x,a})$,
$t\leq s\leq T\}$ to (\ref{BSDEmarkov})--(\ref{cons-markov})
when $(X,I) = \{(X_s^{t,x,a},I_s^{t,a}), t \leq s\leq T\}$.
Moreover, as we shall see in the next paragraph,
this minimal solution is written in this Markovian context as:
$Y_s^{t,x,a} = v(s,X_s^{t,x,a},I_s^{t,x,a})$ where $v$ is the
deterministic function defined on
$[0,T]\times\R^d\times\R^q \rightarrow\R$ by
%
%
\begin{eqnarray}
\label{defv} v(t,x,a) &:=& Y_t^{t,x,a},\qquad (t,x,a)
\in[0,T]\times\R^d\times\R^q.
\end{eqnarray}
We aim at proving that the function $v$ defined by (\ref{defv}) does not
depend actually on its argument $a$, and is a solution in a sense to be
precise to the parabolic IPDE (\ref{IPDE})--(\ref{termIPDE}). Notice that
we do not have a priori any smoothness or even continuity properties on $v$.

To this end, we first recall the definition of (discontinuous)
viscosity solutions to (\ref{IPDE})--(\ref{termIPDE}). For a locally
bounded function $w$
on $[0,T)\times\R^d$, we define its lower semicontinuous (l.s.c. for
short) envelope $w_*$, and upper semicontinuous (u.s.c. for short) envelope
$w^*$ by
\begin{eqnarray*}
w_*(t,x) = \mathop{\operatorname{\lim\inf}\limits
_{(t',x') \rightarrow
(t,x)}}_{t' < T} w
\bigl(t',x'\bigr) \quad\mbox{and}\quad w^*(t,x) =
\mathop{\limsup\limits
_{(t',x') \rightarrow(t,x)}}_{t' < T} %
w\bigl(t',x'
\bigr),
\end{eqnarray*}
for all $(t,x)\in[0,T]\times\R^d$.

%
\begin{Definition}[{[Viscosity solutions to (\ref{IPDE})--(\ref
{termIPDE})]}]\label{devvisocIPDE}
\begin{longlist}[(ii)]
\item[(i)] A function $w$, l.s.c. (resp., u.s.c.) on $[0,T]\times\R^d$,
is called a viscosity supersolution (resp., subsolution) to
(\ref{IPDE})--(\ref{termIPDE}) if
\begin{eqnarray*}
w(T,x) & \geq\mbox{ (resp., $\leq$) }& \sup_{a\in A} g(x,a),
\end{eqnarray*}
for any $x\in\R^d$, and
\begin{eqnarray*}
\biggl( - \frac{\partial\varphi}{\partial t} - \sup_{a\in A} \bigl[
\Lc^a \varphi +f\bigl(\cdot,a,w,\sigma\trans(\cdot,a)D_x
\varphi,\Mc^a\varphi\bigr) \bigr] \biggr) (t,x) & \geq\mbox{ (resp.,
$\leq$) }& 0,
\end{eqnarray*}
for any $(t,x) \in[0,T)\times\R^d$ and any $\varphi\in
C^{1,2}([0,T]\times\R^d)$ such that
\begin{eqnarray*}
(w-\varphi) (t,x) & = & \min_{[0,T]\times\R^d}(w-\varphi)
\qquad \Bigl[\mbox{resp. }\max_{[0,T]\times\R^d}(w-\varphi)\Bigr].
\end{eqnarray*}
\item[(ii)] A locally bounded function $w$ on $[0,T)\times\R^d$ is
called a
viscosity solution to (\ref{IPDE})--(\ref{termIPDE}) if $w_*$ is a
viscosity supersolution
and $w^*$ is a viscosity subsolution to (\ref{IPDE})--(\ref{termIPDE}).
\end{longlist}
\end{Definition}

We can now state the main result of this paper.

%
\begin{Theorem} \label{theomain2}
Assume that conditions \emph{(H$A$)}, \emph{(H$\lambda_\pi$)},
\emph{(HFC)}, \emph{(HBC1)} and \emph{(HBC2)} hold.
The function $v$ in (\ref{defv}) does not depend on the variable $a$ on
$[0,T)\times\R\times\mathring A$, that is,
\begin{eqnarray*}
v(t,x,a) & = & v\bigl(t,x,a'\bigr),\qquad\forall
a,a' \in\mathring A,
\end{eqnarray*}
for all $(t,x)\in[0,T)\times\R^d$. Let us then define by misuse of
notation the function $v$ on $[0,T)\times\R^d$ by
%
%
\begin{eqnarray}
\label{defw} v(t,x) & = & v(t,x,a),\qquad(t,x) \in[0,T)\times\R^d,
\end{eqnarray}
for any $a\in\mathring A$. Then $v$ is a viscosity solution to (\ref
{IPDE}) and a viscosity subsolution to (\ref{termIPDE}). Moreover, if
\emph{(HBC2$'$)} holds, $v$ is a viscosity supersolution to (\ref{termIPDE}).
\end{Theorem}

%
\begin{Remark}
1. Once we have a uniqueness result for the fully nonlinear
IPDE (\ref{IPDE})--(\ref{termIPDE}), Theorem~\ref{theomain2}
provides a
Feynman--Kac representation of this unique solution by means of the
minimal solution to the BSDE
(\ref{BSDEmarkov})--(\ref{cons-markov}). This suggests consequently an
original probabilistic numerical approximation of the nonlinear IPDE
(\ref{IPDE})--(\ref{termIPDE}) by discretization and simulation of
the minimal solution to the BSDE (\ref{BSDEmarkov})--(\ref{cons-markov}).
This issue, especially the treatment of the nonpositive jump constraint,
has been recently investigated in \cite{KLP13a} and \cite{KLP13b},
where the authors analyze the convergence rate of the
approximation scheme, and illustrate their results with several
numerical tests arising for instance in the super-replication of
options in uncertain volatilities and correlations models.
We mention here that a nice feature of our scheme is the fact that the
forward process $(X,I)$ can be easily simulated: indeed, notice that
the jump times of $I$ follow a Poisson distribution of parameter
$\bar\lambda_\pi:= \int_A \lambda_\pi(da)$, and so the pure jump
process $I$ is perfectly simulatable once we know how to simulate the
distribution $\lambda_\pi(da)/\bar\lambda_\pi$ of the jump marks.
Then we can use a standard Euler scheme for simulating the component
$X$. Our scheme does not suffer the curse of dimensionality encountered
in finite difference methods or controlled Markov chains, and takes
advantage of the high dimensional properties of Monte--Carlo methods.

2. We do not address here comparison principles (and so
uniqueness results) for the general parabolic nonlinear IPDE
(\ref{IPDE})--(\ref{termIPDE}).
In the case where the generator function $f(x,a)$ does not depend on
$(y,z,u)$ (see Remark~\ref{remcontrole} below), comparison principle is
proved in \cite{pha98}, and the result can be extended by same
arguments when $f(x,a,y,z)$ depends also on $y,z$ under the Lipschitz
condition {(HBC1)}(ii). When $f$ also depends on $u$, comparison
principle is proved by \cite{barbucpar97} in the semilinear IPDE case,
that is, when $A$ is reduced to a singleton, under condition {(HBC2$'$)}.
We also mention recent results on comparison principles for
IPDE in \cite{barimb08} and references therein.
\end{Remark}

%
\begin{Remark}[(Stochastic control problem)] \label{remcontrole}
1. Let us now consider the particular and important case
where the generator $f(x,a)$ does not depend on $(y,z,u)$. We then observe
that the nonlinear IPDE (\ref{IPDE}) is the Hamilton--Jacobi--Bellman
(HJB) equation associated to the following stochastic control problem:
let us introduce the controlled jump-diffusion process:
%
%
\begin{eqnarray}
\label{controlX} \qquad dX_s^\alpha&=& b\bigl(X_s^\alpha,
\alpha_s\bigr)\,ds + \sigma\bigl(X_s^\alpha,
\alpha_s\bigr)\,dW_s + \int_L
\beta\bigl(X_{s^-}^\alpha,\alpha_s,\ell\bigr)
\tilde\vartheta (ds,d\ell),
\end{eqnarray}
where $W$ is a Brownian motion independent of a random measure
$\vartheta$ on a filtered probability space $(\Omega,\Fc,\F^0,\P)$,
the control $\alpha$ lies in $\Ac_{\F^0}$, the set of $\F
^0$-predictable process valued in $A$, and define the value function
for the control problem
\begin{eqnarray}
w(t,x) := \sup_{\alpha\in\Ac_{\F^0}} \E \biggl[ \int_t^T
f\bigl(X_s^{t,x,\alpha},\alpha_s\bigr)\,ds + g
\bigl(X_T^{t,x,\alpha},\alpha_T\bigr) \biggr],
\nonumber
\\
\eqntext{\displaystyle(t,x) \in[0,T]\times\R^d,}
\end{eqnarray}
where $\{X^{t,x,\alpha}_s,t\leq s\leq T\}$ denotes the solution to
(\ref{controlX}) starting from $x$ at $s = t$, given a control
$\alpha\in\Ac_{\F^0}$. It is well known (see, e.g., \cite{pha98}
or \cite{oksu07})
that the value function $w$ is characterized as the unique
viscosity solution to the dynamic programming HJB equation
(\ref{IPDE})--(\ref{termIPDE}) and, therefore, by Theorem~\ref
{theomain2}, $w = v$.
In other words, we have provided a representation of fully nonlinear
stochastic control problem, including especially control in the
diffusion term, possibly degenerate,
in terms of minimal solution to the BSDE (\ref{BSDEmarkov})--(\ref
{cons-markov}).

2. Combining the BSDE representation of Theorem~\ref
{theomain2} together with the dual representation in Theorem~\ref
{theodual}, we obtain
an original representation for the value function of stochastic control problem
\begin{eqnarray*}
\label{repstodual} &&\sup_{\alpha\in\Ac_{\F^0}} \E \biggl[ \int
_0^T f\bigl(X_t^{\alpha
},
\alpha_t\bigr)\,dt + g\bigl(X_T^{\alpha},
\alpha_T\bigr) \biggr]
\\
&&\qquad= \sup_{\nu\in\Vc_A} \E^\nu \biggl[ \int
_0^T f(X_t,I_t)\, dt+
g(X_T,I_T) \biggr].
\end{eqnarray*}
The right-hand side in the above relation may be viewed as a weak
formulation of
the stochastic control problem. Indeed, it is well known that when
there is only control on the drift, the value function may be
represented in terms of control on change of equivalent probability
measures via Girsanov's theorem for Brownian motion. Such
representation is called weak formulation for stochastic control
problem; see \cite{elk79}.
In the general case, when there is control on the diffusion
coefficient, such ``Brownian'' Girsanov's transformation cannot be
applied, and the idea here is to introduce an exogenous process $I$
valued in the control set $A$, independent of $W$ and $\vartheta$
governing the controlled state process $X^\alpha$, and then to control
the change of equivalent probability measures through a Girsanov's
transformation on this auxiliary process.

3. \emph{Non-Markovian extension.} An interesting issue is to extend
our BSDE representation of stochastic control problem to a
non-Markovian context, that is, when the coefficients $b$, $\sigma$ and
$\beta$ of the controlled process are path-dependent. In this case, we
know from the recent works by Ekren, Touzi and Zhang \cite{ekrtouzha13}
that the value function to the path-dependent stochastic control is a
viscosity solution to a
path-dependent fully nonlinear HJB equation. One possible approach for
getting a BSDE representation to path-dependent stochastic control,
would be
to prove that our minimal solution to the BSDE with nonpositive jumps
is a viscosity solution to the path-dependent fully nonlinear HJB
equation, and then to conclude with a uniqueness result for
path-dependent nonlinear PDE. However, to the best of our knowledge,
there is not yet such comparison result for viscosity supersolution and
subsolution of path-dependent nonlinear PDEs. Instead, we recently
proved in \cite{fuhpha13} by
purely probabilistic arguments that the minimal solution to the BSDE
with nonpositive jumps is equal to the value function of a
path-dependent stochastic control problem, and our approach circumvents
the delicate issue of dynamic programming principle and viscosity
solution in the non-Markovian context. Our result is also obtained
without assuming that $\sigma$ is nondegenerate, in contrast with
\cite{ekrtouzha13}
(see their Assumption~4.7).
\end{Remark}

The rest of this paper is devoted to the proof of Theorem~\ref{theomain2}.

\subsection{Viscosity property of the penalized BSDE}\label{sec3.2}

Let us consider the Markov penalized BSDE associated to (\ref
{BSDEmarkov})--(\ref{cons-markov})
%
%
\begin{eqnarray}
\label{BSDEpen-markov} Y^n_t & = & g(X_T,I_T)
+ \int_t^T f\bigl(X_s,I_s,Y^n_s,Z^n_s,U^n_s
\bigr)\,ds
\nonumber
\\
& & {} + n\int_t^T\int_A
\bigl[R^n_s(a)\bigr]^+\lambda_\pi(da)\,ds- \int_t^T Z^n_s.dW_s
\\
&&{}- \int_t^T \int_L
U^n_s(\ell) \tilde\vartheta(ds,d\ell)-\int_t^T \int_A
R^n_s(a) \tilde\pi (ds,da),
\nonumber
\end{eqnarray}
and denote by $\{
(Y_s^{n,t,x,a},Z_s^{n,t,x,a},U_s^{n,t,x,a},R_s^{n,t,x,a})$, $t\leq
s\leq T\}$ the unique solution to (\ref{BSDEpen-markov}) when
$(X,I) = \{(X_s^{t,x,a},I_s^{t,a}), t \leq s\leq T\}$ for any
initial condition $(t,x,a) \in[0,T]\times\R^d\times\R^q$.
From the Markov
property of the jump-diffusion process $(X,I)$, we recall from \cite
{barbucpar97}
that $Y_s^{n,t,x,a} = v_n(s,X_s^{t,x,a},I_s^{t,a})$,
$t \leq s \leq T$, where $v_n$ is the deterministic function defined on
$[0,T]\times\R^d\times\R^q$ by
%
%
\begin{eqnarray}
\label{defvn} v_n(t,x,a) & := & Y_t^{n,t,x,a},
\qquad(t,x,a) \in[0,T]\times\R ^d\times\R^q.
\end{eqnarray}
From the convergence result (Theorem~\ref{theomain1}) of the penalized
solution, we deduce that the minimal solution of the constrained BSDE is
actually in the form $Y_s^{t,x,a} = v(s,X_s^{t,x,a},I_s^{t,a})$,
$t\leq s\leq T$, with a deterministic function $v$ defined in (\ref{defv}).

Moreover, from the uniform estimate (\ref{bounduni}) and Lemma~\ref
{lemH1}, there exists some positive constant $C$ s.t. for all $n$,
\begin{eqnarray}
&&\bigl\llvert v_n(t,x,a)\bigr\rrvert ^2  \leq C \biggl(
\E\bigl\llvert g\bigl(X_T^{t,x,a},I_T^{t,a}
\bigr)\bigr\rrvert ^2 + \E \biggl[\int_t^T
\bigl\llvert f\bigl(X_s^{t,x,a},I_s^{t,a},0,0,0
\bigr)\bigr\rrvert ^2\,ds \biggr]
\nonumber
\\
&&\hspace*{216.5pt}{}+ \E \Bigl[ \sup_{t\leq s\leq T} \bigl\llvert \bar v
\bigl(s,X_s^{t,x,a}\bigr)\bigr\rrvert ^2 \Bigr]
\biggr),\nonumber
\end{eqnarray}
for all $(t,x,a) \in[0,T]\times\R^d\times\R^q$. From the
polynomial growth condition in {(HBC1)}(i) for $g$ and $f$,
(\ref{polybarv}) for $\bar v$, and the estimate (\ref{estimX}) for
$(X,I)$, we obtain that $v_n$, and thus also $v$ by passing to the
limit, satisfy a polynomial growth condition: there exists some
positive constant $C_v$ and some $p \geq2$, such that for all $n$
%
%
\begin{eqnarray}
\label{polyvvn} \bigl\llvert v_n(t,x,a)\bigr\rrvert + \bigl\llvert
v(t,x,a)\bigr\rrvert & \leq& C_v \bigl(1 + \llvert x\rrvert
^p + \llvert a\rrvert ^p \bigr),
\nonumber
\\[-8pt]
\\[-8pt]
\eqntext{\displaystyle\forall(t,x,a) \in[0,T]\times\R^d\times\R^q.}
\end{eqnarray}

We now consider the parabolic semilinear penalized IPDE for any $n$
%
%
\begin{eqnarray}
\label{IPDEn} \hspace*{10pt}&& {-} \frac{\partial v_n}{\partial t}(t,x,a) - \Lc^a
v_n(t,x,a) - f \bigl(x,a,v_n,\sigma\trans
(x,a)D_x v_n,\Mc^a v_n\bigr)
\\
\hspace*{10pt}&& {}\qquad-\int_A\bigl[v_n
\bigl(t,x,a'\bigr)-v_n(t,x,a)\bigr]\lambda_\pi
\bigl(da'\bigr)
\nonumber
\\
\hspace*{10pt}&&{}\qquad-n\int_A\bigl[v_n
\bigl(t,x,a'\bigr)-v_n(t,x,a)\bigr]^+
\lambda_\pi\bigl(da'\bigr) = 0,
\nonumber
\\
\eqntext{\displaystyle\mbox{on } [0,T)\times\R^d\times\R^q,}
\\
\label{termIPDEn} \hspace*{10pt}&&v_n(T,\cdot,\cdot) = g,\qquad\mbox{on }
\R^d\times\R^q.
\end{eqnarray}

From Theorem~3.4 in Barles et al. \cite{barbucpar97}, we have the
well-known property that the penalized BSDE with jumps (\ref{BSDEpen})
provides a viscosity solution to the penalized IPDE (\ref
{IPDEn})--(\ref{termIPDEn}). Actually, the relation in their paper is
obtained under
{(HBC2$'$)}, which allows the authors to get comparison theorem for
BSDE, but such comparison theorem also holds under the weaker condition
{(HBC2)} as shown in \cite{roy06}, and we then get the following result.

%
\begin{Proposition}
Let Assumptions \emph{(HFC)}, \emph{(HBC1)} and \emph{(HBC2)}
hold. The function $v_n$ in (\ref{defvn}) is a continuous viscosity
solution to
(\ref{IPDEn})--(\ref{termIPDEn}), that is, it is continuous on
$[0,T]\times\R^d\times\R^q$, a viscosity supersolution (resp.,
subsolution) to (\ref{termIPDEn}),
\begin{eqnarray*}
v_n(T,x,a) & \geq\mbox{ (resp., $\leq$) } & g(x,a),
\end{eqnarray*}
for any $(x,a)\in\R^d\times\R^q$, and a viscosity supersolution
(resp., subsolution) to (\ref{IPDEn}),
%
%
\begin{eqnarray}
\label{defviscovn} &&{-} \frac{\partial\varphi}{\partial t}(t,x,a) - \Lc^a
\varphi(t,x,a)
\\
&&{}\qquad-f\bigl(x,a,v_n(t,x,a),\sigma\trans(x,a)D_x
\varphi(t,x,a),\Mc ^a\varphi (t,x,a) \bigr)
\nonumber
\\
&&{}\qquad-\int_A\bigl[\varphi\bigl(t,x,a'
\bigr)-\varphi(t,x,a)\bigr]\lambda_\pi \bigl(da'\bigr)
\nonumber
\\
&&{}\qquad-n\int_A\bigl[\varphi\bigl(t,x,a'
\bigr)-\varphi(t,x,a)\bigr]^+\lambda_\pi\bigl(da'\bigr)
\geq\mbox{ (resp., $\leq$) }0,
\nonumber
\end{eqnarray}
for any $(t,x,a) \in[0,T)\times\R^d\times\R^q$ and any
$\varphi
\in C^{1,2}([0,T]\times(\R^d\times\R^q))$ such that
%
%
\begin{eqnarray}
\label{phiglo} 
(v_n-\varphi) (t,x,a) = \min
_{[0,T]\times\R^d\times\R
^q}(v_n-\varphi)
\nonumber
\\[-8pt]
\\[-8pt]
\eqntext{\displaystyle\Bigl[\mbox{resp., }\max_{[0,T]\times\R^d\times
\R^q}(v_n-
\varphi) \Bigr].}
\end{eqnarray}
\end{Proposition}

In contrast to local PDEs with no integro-differential terms, we cannot
restrict in general the global minimum (resp., maximum) condition on
the test functions for the definition of viscosity supersolution
(resp., subsolution) to local minimum (resp., maximum) condition. In
our IPDE case,
the nonlocal terms appearing in (\ref{IPDEn}) involve the values w.r.t.
the variable $a$ only on the set $A$. Therefore, we are able to
restrict the
global extremum condition on the test functions to extremum on
$[0,T]\times\R^d\times A$. More precisely, we have the following
equivalent definition of viscosity solutions, which will be used later.

%
\begin{Lemma} \label{lemviscovnA}
Assume that \emph{(H$\lambda_\pi$)}, \emph{(HFC)}, and
\emph{(HBC1)} hold. In the definition of $v_n$ being a viscosity
supersolution (resp., subsolution) to (\ref{IPDEn}) at a point
$(t,x,a) \in[0,T)\times\R^d\times\mathring A$, we can replace
condition (\ref{phiglo}) by
\begin{eqnarray*}
0 = (v_n-\varphi) (t,x,a) & = & \min_{[0,T]\times\R^d\times
\mathring A}(v_n-
\varphi)\qquad \Bigl[\mbox{resp., }\max_{[0,T]\times\R^d\times\mathring A}(v_n-
\varphi ) \Bigr],
\end{eqnarray*}
and suppose that the test function $\varphi$ is in
$C^{1,2,0}([0,T]\times\R^d\times\R^q)$.
\end{Lemma}

\begin{pf} We treat only the supersolution case as the subsolution
case is proved by same arguments, and proceed in two steps.

\noindent\emph{Step 1.} Fix $(t,x,a) \in[0,T)\times\R^d\times\R
^q$, and
let us show that the viscosity supersolution inequality (\ref{defviscovn})
also holds for any
test function $\varphi$ in $C^{1,2,0}([0,T]\times\R^d\times\R^q)$ s.t.
%
%
\begin{eqnarray}
\label{condminglob1} 
(v_n-\varphi) (t,x,a) & = & \min
_{[0,T]\times\R^d\times\R
^q}(v_n-\varphi ).
\end{eqnarray}
We may assume w.l.o.g. that the minimum for such test function $\varphi
$ is zero, and let us define for $r>0$ the function $\varphi^r$ by
%
\begin{eqnarray*}
\varphi^r\bigl(t',x',a'
\bigr) & = & \varphi\bigl(t',x',a'\bigr)
\biggl(1-\Phi \biggl(\frac
{\llvert x'\rrvert ^{2}+\llvert a'\rrvert ^{2}}{r^2} \biggr) \biggr)
\\
&&{} - C_v \Phi \biggl(\frac{\llvert x'\rrvert ^{2}+\llvert a'\rrvert ^{2}}{r^2} \biggr) \bigl(1+\bigl
\llvert x'\bigr\rrvert ^{p}+\bigl\llvert a'
\bigr\rrvert ^{p} \bigr),
\end{eqnarray*}
where $C_v > 0$ and $p \geq2$ are the constant and degree
appearing in the polynomial growth condition (\ref{polyvvn}) for $v_n$,
$\Phi\dvtx\R_+\rightarrow[0,1]$ is a function in $C^\infty(\R_+)$
such that
$\Phi|_{[0,1]}\equiv0$ and $\Phi|_{[2,+\infty)}\equiv1$.
Notice that $\varphi^r\in C^{1,2,0}([0,T]\times\R^d\times\R^q)$,
%
%
\begin{eqnarray}
\label{convphikphir} \bigl(\varphi^r,D_x\varphi^r,D^2_x
\varphi^r\bigr) & \longrightarrow& \bigl(\varphi,D_x
\varphi,D_x^2\varphi\bigr)\qquad\mbox{as } r\rightarrow
\infty
\end{eqnarray}
locally uniformly on $[0,T]\times\R^d\times\R^q$, and that there exists
a constant $C_r>0$ such that
%
\begin{eqnarray}
\label{croisslinfcttestr} \bigl\llvert \varphi^r\bigl(t',x',a'
\bigr)\bigr\rrvert & \leq& C_r \bigl(1+\bigl\llvert x'
\bigr\rrvert ^{p}+\bigl\llvert a'\bigr\rrvert
^{p} \bigr)
\end{eqnarray}
for all $(t',x',a')\in[0,T]\times\R^q\times\R^d$. Since $\Phi$ is
valued in $[0,1]$, we deduce from the polynomial growth condition (\ref
{polyvvn}) satisfied by $v_n$ and (\ref{condminglob1}) that
$\varphi^r \leq v_n$ on $[0,T]\times\R^d\times\R^q$
for all $r> 0$.
Moreover, we have $\varphi^r(t,x,a) = \varphi(t,x,a)$ [$=
v_n(t,x,a)$] for $r$ large enough. Therefore, we get 
%
\begin{eqnarray}
\label{condstricminvn} 
\bigl(v_n-\varphi^r\bigr)
(t,x,a) & = & \min_{[0,T]\times\R^d\times\R
^q}\bigl(v_n-
\varphi^r\bigr),
\end{eqnarray}
for $r$ large enough,
and we may assume w.l.o.g. that this minimum is strict.
Let $(\varphi_k^{r})_k$ be a sequence of function in
$C^{1,2}([0,T]\times(\R^d\times\R^q))$ satisfying (\ref
{croisslinfcttestr}) and such that
%
%
\begin{eqnarray}
\label{conv1convr} \bigl(\varphi_k^{r},D_x
\varphi_k^{r}, D^2_x
\varphi_k^{r}\bigr) & \longrightarrow & \bigl(
\varphi^{r},D_x\varphi^{r},
D^2_x\varphi^{r}\bigr)\qquad\mbox{as } k
\rightarrow\infty,
\end{eqnarray}
locally uniformly on $[0,T]\times\R^d\times\R^q$. From the growth
conditions (\ref{polyvvn}) and (\ref{croisslinfcttestr}) on the
continuous functions
$v_n$ and $\varphi_k^r$, we can assume w.l.o.g. [up to an usual
negative perturbation of the function $\varphi_r^k$ for large $(x',a')$],
that there exists a bounded sequence $(t_k,x_k,a_k)_k$ in $[0,T]\times
\R
^d\times\R^q$ such that
%
\begin{eqnarray}
\label{minphikr1} 
\bigl(v_n-\varphi^{r}_k
\bigr) (t_k,x_k,a_k) & = & \min
_{[0,T]\times\R^d\times
\R
^q}\bigl(v_n-\varphi^{r}_k
\bigr).
\end{eqnarray}
The sequence $(t_k,x_k,a_k)_k$ converges up to a subsequence, and thus,
by (\ref{condstricminvn}), (\ref{conv1convr}) and (\ref{minphikr1}),
we have\vspace*{1pt}
%
%
\begin{eqnarray}
\label{convtxar} (t_k,x_k,a_k) &
\rightarrow& (t,x,a),\qquad\mbox{as } k\rightarrow\infty.
\end{eqnarray}
Now, from the viscosity supersolution property of $v_n$ at
$(t_k,x_k,a_k)$ with the test function $\varphi_k^{r}$, we
have
\begin{eqnarray*}
&&{-} \frac{\partial\varphi_k^{r}}{\partial t}(t_k,x_k,a_k) -
\Lc ^{a_k} \varphi _k^{r}(t_k,x_k,a_k)
\\
&&{}\qquad-f\bigl(x_k,a_k,v_n(t_k,x_k,a_k),
\sigma\trans (x_k,a_k)D_x\varphi
_k^{r}(t_k,x_k,a_k),
\Mc^{a_k}\varphi_k^{r}(t_k,x_k,a_k)
\bigr)
\\
&&{}\qquad-\int_A\bigl[\varphi_k^{r}
\bigl(t_k,x_k,a'\bigr)-\varphi
_k^{r}(t_k,x_k,a_k)
\bigr]\lambda _\pi\bigl(da'\bigr)
\\
&&{}\qquad-n\int_A\bigl[\varphi_k^{r}
\bigl(t_k,x_k,a'\bigr)-\varphi
_k^{r}(t_k,x_k,a_k)
\bigr]^+\lambda _\pi\bigl(da'\bigr) \geq0.
\end{eqnarray*}
Sending $k$ and $r$ to infinity, and using (\ref{convphikphir}), (\ref
{conv1convr}) and (\ref{convtxar}), we obtain the viscosity
supersolution inequality at $(t,x,a)$ with the test function $\varphi$.

\noindent\textit{Step 2.} Fix $(t,x,a) \in[0,T)\times\R^d\times
\mathring
A$, and let $\varphi$ be a test function in
$C^{1,2}([0,T]\times(\R^d\times\R^q))$ such that
%
%
\begin{eqnarray}
\label{phiglo2} 0 = (v_n-\varphi) (t,x,a) & = & \min
_{[0,T]\times\R^d\times
\mathring
A}(v_n-\varphi).
\end{eqnarray}
By the same arguments as in (\ref{croisslinfcttestr}), we can assume
w.l.o.g. that $\varphi$ satisfies the polynomial growth condition
\begin{eqnarray*}
\bigl\llvert \varphi\bigl(t',x',a'
\bigr)\bigr\rrvert & \leq& C\bigl(1+\bigl\llvert x'\bigr\rrvert
^{p}+\bigl\llvert a'\bigr\rrvert ^{p}\bigr),
\qquad\bigl(t',x',a'\bigr) \in[0,T]
\times\R^d\times\R^q,
\end{eqnarray*}
for some positive constant $C$. Together with (\ref{polyvvn}), and since
$A$ is compact, we have
%
%
\begin{eqnarray}
\label{maju-phi} (v_n-\varphi) \bigl(t',x',a'
\bigr) & \geq& - C\bigl(1+\bigl\llvert x'\bigr\rrvert
^{p} +\bigl\llvert d_{A}\bigl(a'\bigr)\bigr
\rrvert ^{p}\bigr),
\end{eqnarray}
for all $(t',x',a')\in[0,T]\times\R^d\times\R^q$, where $d_A(a')$ is
the distance from $a'$ to $A$.
Fix $\eps>0$ and
define the function $ \varphi_\eps\in C^{1,2,0}([0,T]\times\R
^d\times\R
^q)$ by
\begin{eqnarray*}
\varphi_\eps\bigl(t',x',a'
\bigr) & = & \varphi\bigl(t',x',a'\bigr)
- \Phi \biggl(\frac
{d_{A_{\eps}}(a')}{{\eps}} \biggr)C\bigl(1+\bigl\llvert x'
\bigr\rrvert ^{p}+\bigl\llvert d_{A}\bigl(a'
\bigr)\bigr\rrvert ^{p}\bigr)
\end{eqnarray*}
for all $(t',x',a')\in[0,T]\times\R^d\times\R^q$, where
%
%
\begin{eqnarray}
\label{defAeps} A_{\eps} & = & \bigl\{a'\in A \dvtx
d_{\partial A}\bigl(a'\bigr) \geq{\eps } \bigr\},
\end{eqnarray}
and $\Phi\dvtx\R_+\rightarrow[0,1]$ is
a function in $C^\infty(\R_+)$ such that $\Phi|_{[0,{1/2}]}\equiv
0$ and $\Phi|_{[1,+\infty)}\equiv1$. Notice that
%
%
\begin{eqnarray}
\label{convphir} \bigl(\varphi_\eps,D_x\varphi_\eps,D^2_x
\varphi_\eps\bigr) & \longrightarrow& \bigl(\varphi,D_x
\varphi,D_x^2\varphi\bigr) \qquad\mbox{as } \eps
\rightarrow0,
\end{eqnarray}
locally uniformly on $[0,T]\times\R^d\times\mathring A$. We notice from
(\ref{maju-phi}) and the definition of $\varphi_\eps$ that
$\varphi_\eps\leq v_n$ on $[0,T]\times\R^d\times A_{\eps}^c$.
Moreover, since $\varphi_\eps\leq\varphi$ on $[0,T]\times\R
^d\times\R
^q$, $\varphi_\eps=\varphi$ on $[0,T]\times\R^d\times\mathring
A_{\eps
}$ and $a\in\mathring A$, we get by (\ref{phiglo2}) for $\eps$ small enough
%
\begin{eqnarray*}
0 = (v_n-\varphi_\eps) (t,x,a) & = & \min
_{[0,T]\times\R^d\times\R
^q}(v_n-\varphi_\eps).
\end{eqnarray*}
From Step 1, we then have
\begin{eqnarray*}
&&{-} \frac{\partial\varphi_\eps}{\partial t}(t,x,a) - \Lc ^a\varphi_\eps(t,x,a)
\nonumber
\\
&&{}\qquad-f\bigl(x,a,v_n(t,x,a),\sigma\trans(x,a)D_x
\varphi_\eps (t,x,a),\Mc ^a\varphi _\eps(t,x,a)
\bigr)
\\
&&{}\qquad-\int_A\bigl[\varphi_\eps
\bigl(t,x,a'\bigr)-\varphi_\eps (t,x,a)\bigr]
\lambda_\pi\bigl(da'\bigr)
\\
&&{}\qquad -n\int_A\bigl[\varphi_\eps
\bigl(t,x,a'\bigr)-\varphi_\eps(t,x,a)\bigr]^+
\lambda_\pi\bigl(da'\bigr) \geq0 .
\end{eqnarray*}
By sending $\eps$ to zero with (\ref{convphir}), and using $a\in
\mathring A$ with {(H$\lambda_\pi$)}(ii), we get the required
viscosity
subsolution inequality at $(t,x,a)$ for the test function $\varphi$.
\end{pf}

\subsection{The nondependence of the function $v$ in the variable
$a$}\label{secnondepa}

In this subsection, we aim to prove that the function $v(t,x,a)$ does
not depend on $a$.
From the relation defining the Markov BSDE (\ref{BSDEmarkov}), and since
for the minimal solution
$(Y^{t,x,a},Z^{t,x,a},U^{t,x,a},R^{t,x,a},K^{t,x,a})$ to (\ref
{BSDEmarkov})--(\ref{cons-markov}), the process $K^{t,x,a}$ is predictable,
we observe that the $A$-jump component $R^{t,x,a}$ is expressed in
terms of $Y^{t,x,a} = v(\cdot,X^{t,x,a},I^{t,x,a})$ as
\begin{eqnarray*}
R_s^{t,x,a}\bigl(a'\bigr) &=& v
\bigl(s,X_{s^-}^{t,x,a},a'\bigr) - v
\bigl(s,X_{s^-}^{t,x,a},I_{s^-}^{t,x,a}\bigr),
\qquad t \leq s \leq T, a ' \in A,
\end{eqnarray*}
for all $(t,x,a) \in[0,T]\times\R^d\times\R^q$. From the
$A$-nonpositive constraint (\ref{cons-markov}), this yields
\begin{eqnarray*}
\E \biggl[ \int_t^{t+h} \int
_A \bigl[ v\bigl(s,X_{s}^{t,x,a},a'
\bigr) - v\bigl(s,X_{s}^{t,x,a},I_{s}^{t,x,a}
\bigr) \bigr]^+ \lambda_\pi\bigl(da'\bigr)\,ds \biggr]
&=& 0,
\end{eqnarray*}
for any $h > 0$.
If we knew a priori that the function $v$ was continuous on
$[0,T)\times
\R^d\times A$, we could obtain by sending $h$ to zero in the above
equality divided by $h$ (and by dominated convergence theorem), and
from the mean-value theorem
\begin{eqnarray*}
\int_A \bigl[ v\bigl(t,x,a'\bigr)-v(t,x,a)
\bigr]^+ \lambda_\pi\bigl(da'\bigr) &=& 0.
\end{eqnarray*}
Under condition {(H$\lambda_\pi$)}(i), this would prove that
$v(t,x,a)\geq v(t,x,a')$ for any $a,a'\in A$, and thus the function $v$
would not depend on $a$ in $A$.

Unfortunately, we are not able to prove directly the continuity of $v$
from its very definition (\ref{defv}),
and instead, we shall rely on viscosity solutions approach to derive
the nondependence of $v(t,x,a)$ in $a \in\mathring A$. To this
end, let us introduce the following first-order PDE:
%
%
\begin{eqnarray}
\label{PDEnondepva} - \bigl\llvert D_av(t,x,a)\bigr\rrvert & = & 0 ,
\qquad(t,x,a)\in[0,T)\times\R ^d\times \mathring A.
\end{eqnarray}

%
\begin{Lemma}\label{lemderparanulsolviscoglob}
Let assumptions \emph{(H$\lambda_\pi$)}, \emph{(HFC)}, \emph{(HBC1)}
and \emph{(HBC2)} hold.
The function $v$ is a viscosity supersolution to (\ref{PDEnondepva}):
for any $(t,x,a)\in[0,T)\times\R^d\times\mathring{A}$ and any function
$\varphi\in C^{1,2}([0,T]\times(\R^d\times\R^q))$ such that
$(v-\varphi
)(t,x,a) = \min_{[0,T]\times\R^d\times\R^q}(v-\varphi)$,
we have
\begin{eqnarray*}
- \bigl\llvert D_a \varphi(t,x,a)\bigr\rrvert &\geq& 0, \quad
\mbox{that is} \quad D_a \varphi (t,x,a) = 0.
\end{eqnarray*}
\end{Lemma}

\begin{pf} We know that $v$ is the pointwise limit of the
nondecreasing sequence of functions $(v_n)$. By continuity of $v_n$,
the function $v$ is l.s.c. and we have (see, e.g., \cite{bar94}, page 91)
%
%
\begin{eqnarray}
\label{vinf1} v = v_* & = & \operatorname{\lim\inf}\limits
_{n\rightarrow\infty
}{}_* v_n,
\end{eqnarray}
where
\begin{eqnarray}
\liminf_{n\rightarrow\infty} {}_*v_n(t,x,a) := \mathop{\mathop {
\liminf_{n\rightarrow\infty}}\limits
_{(t',x',a')\rightarrow
(t,x,a)}}_{t'<T} %
v_n\bigl(t',x',a'\bigr),
\nonumber
\\
\eqntext{\displaystyle(t,x,a)\in[0,T]\times\R^d\times
\R^q.}
\end{eqnarray}
Let $(t,x,a) \in[0,T)\times\R^d\times\mathring{A}$, and
$\varphi
\in C^{1,2}([0,T]\times(\R^d\times\R^q))$, such that
$(v-\varphi)(t,x,a) = \min_{[0,T]\times\R^d\times\R
^q}(v-\varphi
)$. We may assume w.l.o.g. that this minimum is strict:
%
\begin{eqnarray}
\label{strictmin1} (v-\varphi) (t,x,a) & = & \mbox{strict}\min_{[0,T]\times\R^d\times
\R
^q}(v-
\varphi).
\end{eqnarray}
Up to a suitable negative perturbation of $\varphi$ for large $(x,a)$,
we can assume w.l.o.g. that there exists a bounded sequence
$(t_n,x_n,a_n)_n$ in
$[0,T]\times\R^d\times\R^q$ such that
%
%
\begin{eqnarray}
\label{minvnmoinsvarphialpha1} (v_n-\varphi) (t_n,x_n,a_n)
& = & \min_{[0,T]\times\R^d\times\R
^q}(v_n-\varphi).
\end{eqnarray}
From (\ref{vinf1}), (\ref{strictmin1}) and (\ref{minvnmoinsvarphialpha1}),
we then have, up to a subsequence
%
%
\begin{eqnarray}
\label{conv suiteoptimisante1} \bigl(t_n,x_n,a_n,v_n(t_n,x_n,a_n)
\bigr) & \longrightarrow& \bigl(t,x,a,v(t,x,a)\bigr)\qquad\mbox{as } n \rightarrow
\infty.
\end{eqnarray}
%
Now, from the viscosity supersolution property of $v_n$ at
$(t_n,x_n,a_n)$ with the test function $\varphi$, we have by
(\ref{minvnmoinsvarphialpha1}),
\begin{eqnarray*}
&&{-} \frac{\partial\varphi}{\partial t}(t_n,x_n,a_n) - \Lc
^{a_n}\varphi(t_n,x_n,a_n)
\\
&&{}\qquad-f\bigl(x_n,a_n,v_n(t_n,x_n,a_n),
\sigma\trans (x_n,a_n)D_x\varphi
(t_n,x_n,a_n),\Mc^{a_n}
\varphi(t_n,x_n,a_n)\bigr)
\\
&&{}\qquad-\int_A\bigl[\varphi\bigl(t_n,x_n,a'
\bigr)-\varphi(t_n,x_n,a_n)\bigr]\lambda
_\pi\bigl(da'\bigr)
\\
&&{}\qquad-n\int_A\bigl[\varphi\bigl(t_n,x_n,a'
\bigr)-\varphi (t_n,x_n,a_n)\bigr]^+
\lambda_\pi \bigl(da'\bigr) \geq0,
\end{eqnarray*}
which implies
\begin{eqnarray*}
& & \int_A\bigl[\varphi\bigl(t_n,x_n,a'
\bigr)-\varphi(t_n,x_n,a_n)\bigr]^+
\lambda_\pi\bigl(da'\bigr)
\\
&&\qquad\leq{1\over n} \biggl[ - \frac{\partial\varphi}{\partial
t}(t_n,x_n,a_n)
- \Lc ^{a_n}\varphi (t_n,x_n,a_n)
\\
& &\hspace*{45pt} {} - f\bigl(x_n,a_n,v_n(t_n,x_n,a_n),
\sigma\trans (x_n,a_n) D_x\varphi(t_n,x_n,a_n),
\\
&&\hspace*{207pt} {}
\Mc^{a_n}\varphi (t_n,x_n,a_n)\bigr)
\\
& &\hspace*{98pt} {} - \int_A\bigl[\varphi
\bigl(t_n,x_n,a'\bigr)-\varphi
(t_n,x_n,a_n)\bigr]\lambda_\pi
\bigl(da'\bigr) \biggr].
\end{eqnarray*}
Sending $n$ to infinity, we get from (\ref{conv suiteoptimisante1}), the
continuity of coefficients $b,\sigma,\beta$ and $f$, and the dominated
convergence theorem
\begin{eqnarray*}
\int_A \bigl[\varphi\bigl(t,x,a'\bigr)-
\varphi(t,x,a)\bigr]^+\lambda_\pi\bigl(da'\bigr) & = & 0
.
\end{eqnarray*}
%
Under {(H$\lambda_\pi$)}, this means that $\varphi
(t,x,a)=\max_{a'\in A}\varphi(t,x,a')$. Since $a\in\mathring A$, we deduce that
$D_a\varphi(t,x,a) = 0$.
\end{pf}

We notice that the PDE (\ref{PDEnondepva}) involves only differential
terms in the variable $a$. Therefore, we can freeze the terms $(t,x)\in
[0,T)\times\R^d$ in the PDE (\ref{PDEnondepva}), that is, we can take
test functions not depending on the variables $(t,x)$ in the definition
of the viscosity solution, as shown in the following lemma.

%
\begin{Lemma}\label{lemisolevarasursol}
Let assumptions \emph{(H$\lambda_\pi$)}, \emph{(HFC)}, \emph{(HBC1)}
and \emph{(HBC2)} hold.
For any $(t,x)\in[0,T)\times\R^d$, the function $v(t,x,\cdot)$ is a
viscosity supersolution to
\begin{eqnarray*}
- \bigl\llvert D_av(t,x,a)\bigr\rrvert & = & 0 ,\qquad a\in
\mathring A,
\end{eqnarray*}
that is, for any $a\in\mathring{A}$ and any function $\varphi\in
C^{2}(\R^q)$ such that
$(v(t,x,\cdot)-\varphi)(a) = \min_{\R^q}(v(t,x,\cdot)-\varphi)$,
we have:
$- |D_a \varphi(a) | \geq0$ (and so $= 0$).
\end{Lemma}

\begin{pf}
Fix $(t,x)\in[0,T)\times\R^d$, $a\in\mathring{A}$
and $\varphi\in C^2(\R^q)$ such that
%
%
\begin{eqnarray}
\label{condminpartiellem2} \bigl(v(t,x,\cdot)-\varphi\bigr) (a) & = & \min
_{\R^q}\bigl(v(t,x,\cdot)-\varphi\bigr).
\end{eqnarray}
As usual, we may assume w.l.o.g. that this minimum is strict and that
$\varphi$ satisfies the growth condition $\sup_{a'\in\R^q}{|\varphi
(a')|\over1+|a'|^p}<\infty$.
Let us then define for $n\geq1$, the function $\varphi^n\in
C^{1,2}([0,T]\times(\R^d\times\R^q))$ by
\begin{eqnarray*}
\varphi^n\bigl(t',x',a'
\bigr) & = & \varphi\bigl(a'\bigr)-n \bigl(\bigl\llvert
t'-t\bigr\rrvert ^2+\bigl\llvert x'-x
\bigr\rrvert ^{2p} \bigr)- \bigl\llvert a'-a\bigr\rrvert
^{2p}
\end{eqnarray*}
for all $(t',x',a')\in[0,T]\times\R^d\times\R^q$. From the growth
condition (\ref{polyvvn}) on the l.s.c. function $v$, and the growth condition
on the continuous function $\varphi$, one can find for any $n\geq1$ an
element $(t_n,x_n,a_n)$ of $[0,T]\times\R^d\times\R^q$ such that
\begin{eqnarray*}
\bigl(v-\varphi^n\bigr) (t_n,x_n,a_n)
& = & \min_{[0,T]\times\R^d\times\R
^q}\bigl(v-\varphi^n\bigr).
\end{eqnarray*}
In particular, we have
%
%
\begin{eqnarray}
\label{condminlem2} v(t,x,a)-\varphi(a) & = & \bigl(v-\varphi^n\bigr)
(t,x,a) \geq\bigl(v-\varphi ^n\bigr) (t_n,x_n,a_n)
\nonumber
\\
& = & v(t_n,x_n,a_n) - \varphi(a_n)
\nonumber
\\
&&{} + n\bigl(\llvert t_n-t\rrvert ^2+\llvert
x_n-x\rrvert ^{2p}\bigr)+\llvert a_n-a\rrvert
^{2p}
\\
& \geq& v(t_n,x_n,a_n) -
v(t,x,a_n) + v(t,x,a)-\varphi(a)
\nonumber
\\
& &{} + n\bigl(\llvert t_n-t\rrvert ^2+\llvert
x_n-x\rrvert ^{2p}\bigr)+\llvert a_n-a\rrvert
^{2p}
\nonumber
\end{eqnarray}
by (\ref{condminpartiellem2}),
which implies from the growth condition (\ref{polyvvn}) on $v$
\begin{eqnarray*}
n\bigl(\llvert t_n-t\rrvert ^2+\llvert x_n-x
\rrvert ^{2p}\bigr)+\llvert a_n-a\rrvert ^{2p} &
\leq& C\bigl(1+\llvert x_n-x\rrvert ^{p}+\llvert
a_n-a\rrvert ^p\bigr).
\end{eqnarray*}
Therefore, the sequences $\{n(\llvert t_n-t\rrvert ^2+\llvert x_n-x\rrvert ^{2p})\}_n$ and
$(\llvert a-a_n\rrvert ^{2p})_n$ are bounded and (up to a subsequence) we have:
$(t_n,x_n,a_n) \longrightarrow(t,x,a_\infty)$ as $n$ goes to infinity,
for some $a_\infty\in\R^q$. Actually, since $v(t,x,a)-\varphi(a)
\geq
v(t_n,x_n,a_n)-\varphi(a_n)$ by (\ref{condminlem2}), we obtain by
sending $n$ to infinity and since the minimum in (\ref
{condminpartiellem2}) is strict, that $a_\infty= a$, and so
\begin{eqnarray*}
(t_n,x_n,a_n) & \longrightarrow& (t,x,a)
\qquad\mbox{as } n\rightarrow \infty.
\end{eqnarray*}
On the other hand, from Lemma~\ref{lemderparanulsolviscoglob} applied
to $(t_n,x_n,a_n)$ with the test function $\varphi^n$, we have
\begin{eqnarray*}
0 = D_a\varphi^n(t_n,x_n,a_n)
&=& D_a\varphi(a_n)- 2p (a_n-a)\llvert
a_n-a\rrvert ^{2p-1},
\end{eqnarray*}
for all $n\geq1$. Sending $n$ to infinity we get the required result:
$D_a\varphi(a) = 0$.
\end{pf}

We are now able to state the main result of this subsection.

%
\begin{Proposition}\label{THMnondepva}
Let assumptions \emph{(H$A$)}, \emph{(H$\lambda_\pi$)}, \emph{(HFC)},
\emph{(HBC1)} and \emph{(HBC2)} hold.
The function $v$ does not depend on the variable $a$ on $[0,T)\times\R
^d \times\mathring{A}$:
\begin{eqnarray*}
v(t,x,a) & = & v\bigl(t,x,a'\bigr),\qquad a,a' \in
\mathring A,
\end{eqnarray*}
for any $(t,x)\in[0,T)\times\R^d$.
\end{Proposition}

\begin{pf} We proceed in four steps.

\noindent\textit{Step 1.} \textit{Approximation by inf-convolution.}
We introduce the family of functions $(u_n)_n$ defined by
\begin{eqnarray*}
u_n(t,x,a) & = & \inf_{a'\in A} \bigl[ v
\bigl(t,x,a'\bigr)+{ n}\bigl\llvert a-a'\bigr\rrvert
^{2p} \bigr],\qquad (t,x,a) \in[0,T]\times\R^d\times A.
\end{eqnarray*}
It is clear that the sequence $(u_n)_n$ is nondecreasing and
upper-bounded by $v$. Moreover, since $v$ is l.s.c., 
we have the pointwise convergence of $u_n$ to $v$ on $[0,T]\times\R
^d\times A$.
Indeed, fix some $(t,x,a)\in[0,T]\times\R^d\times A$. Since $v$ is l.s.c.,
there exists a sequence $(a_n)_n$ valued in $A$ such that
\begin{eqnarray*}
u_n(t,x,a) & = & v(t,x,a_n)+{ n}\llvert
a-a_n\rrvert ^{2p},
\end{eqnarray*}
for all $n\geq1$.
Since $A$ is compact, the sequence $(a_n)$ converges, up to a
subsequence, to some $a_\infty\in A$. Moreover,
since $u_n$ is upper-bounded by $v$ and $v$ is l.s.c., we see that
$a_\infty= a$ and
%
%
\begin{eqnarray}
\label{convun} u_n(t,x,a) & \longrightarrow& v(t,x,a)\qquad\mbox{as }
n\rightarrow \infty.
\end{eqnarray}

\noindent\textit{Step 2.} \textit{A test function for $u_n$ seen as
a test
function for $v$.}
For $r,\delta>0$ let us define the integer $N(r,\delta)$ by
\begin{eqnarray*}
N(r,\delta) & = & \min \biggl\{ n\in\N\dvtx n \geq\frac{2C_v(1+
2^{2p-5}+ r^p+2^{p-1}\max_{a\in A}\llvert a\rrvert ^p)}{ ({\delta}/{2}
)^{2p}}
+C_v \biggr\},
\end{eqnarray*}
where $C_v$ is the constant in the growth condition (\ref{polyvvn}), and
define the set $\mathring{A}_\delta$ by
\begin{eqnarray*}
\mathring{A}_\delta& = & \Bigl\{a\in A \dvtx d(a,\partial A):=\min
_{a'\in
\partial A}\bigl\llvert a-a'\bigr\rrvert >\delta \Bigr
\}.
\end{eqnarray*}
Fix $(t,x) \in[0,T)\times\R^d$. We now prove that for any
$\delta
>0$, $n\geq N(\llvert x\rrvert,\delta)$, $a\in\mathring{A}_\delta$ and
$\varphi\in C^{2}(\R^q)$ such that
%
%
\begin{eqnarray}
\label{min0un-phi} 0 = \bigl(u_n(t,x,\cdot)-\varphi\bigr) (a) & = & \min
_{\R^q}\bigl(u_n(t,x,\cdot )-\varphi\bigr),
\end{eqnarray}
there exists $a_n\in\mathring{A}$ and $\psi\in C^2(\R^q)$
such that
%
%
\begin{eqnarray}
\label{min0upsi} 0 = \bigl(v(t,x,\cdot)-\psi\bigr) (a_n) & = & \min
_{\R^q}\bigl(v(t,x,\cdot)-\psi\bigr),
\end{eqnarray}
and
%
%
\begin{eqnarray}
\label{derivpsi} D_a \psi(a_n) &=& D_a
\varphi(a).
\end{eqnarray}
To this end, we proceed in two substeps.

\noindent\textit{Substep 2.1.} We prove that for any $\delta>0$,
$(t,x,a)\in
[0,T)\times\R^d\times\mathring{A}_\delta$, and any $n\geq
N(\llvert x\rrvert,\delta)$
\begin{eqnarray*}
\operatorname{\arg\min}\limits
_{a'\in A} \bigl\{ v\bigl(t,x,a'\bigr)+n
\bigl\llvert a'-a\bigr\rrvert ^{2p} \bigr\} & \subset&
\mathring{A}.
\end{eqnarray*}
Fix $(t,x,a)\in[0,T)\times\R^d\times\mathring{A}_\delta$ and let
$a_n\in A$ such that
\begin{eqnarray*}
v(t,x,a_n)+n\llvert a_n-a\rrvert ^{2p} & = &
\min_{a'\in A} \bigl[ v\bigl(t,x,a'\bigr)+n\bigl
\llvert a'-a\bigr\rrvert ^{2p} \bigr] .
\end{eqnarray*}
Then we have
\begin{eqnarray*}
v(t,x,a_n)+n\llvert a_n-a\rrvert ^{2p} & \leq&
v(t,x,a),
\end{eqnarray*}
and by (\ref{polyvvn}), this gives
\begin{eqnarray*}
&&{-} C_v\Bigl(1+\llvert x\rrvert ^p+ 2^{p-1}
\max_{a\in A} \llvert a\rrvert ^p + 2^{p-1}
\llvert a_n-a\rrvert ^p\Bigr)+n\llvert a_n-a
\rrvert ^{2p}
\\
&& \qquad\leq C_v\bigl(1+\llvert x\rrvert ^p+\llvert a
\rrvert ^p\bigr).
\end{eqnarray*}
Then using the inequality $2\alpha\beta\leq\alpha^2+\beta^2$ to
the product $2\alpha\beta= 2^{p-1}\llvert a_n-a\rrvert ^p$, we get
\begin{eqnarray*}
(n-C_v)\llvert a_n-a\rrvert ^{2p} & \leq&
2C_v\Bigl(1+ 2^{2p-5}+ \llvert x\rrvert ^p+2^{p-1}
\max_{a\in
A}\llvert a\rrvert ^p\Bigr).
\end{eqnarray*}
For $n\geq N(\llvert x\rrvert,\delta)$, we get from the definition
of $N(r,\delta)$
\begin{eqnarray*}
\llvert a_n-a\rrvert & \leq& {\delta\over2},
\end{eqnarray*}
which shows that $a_n\in\mathring{A}$ since $a \in\mathring
{A}_\delta$.

\noindent\textit{Substep 2.2.} Fix $\delta>0$, $(t,x,a)\in
[0,T)\times\R^d\times
\mathring{A}_\delta$, and $\varphi\in C^{2}(\R^q)$ satisfying (\ref
{min0un-phi}).
Let us then choose $a_n\in\arg\min\{v(t,x,a')+n\llvert a'-a\rrvert ^{2p}\dvtx
a'\in
A\}
$, and define $\psi\in C^{2}(\R^q)$ by
\begin{eqnarray*}
\psi\bigl(a'\bigr) & = & \varphi\bigl(a+ a'-a_n
\bigr)-n\llvert a_n-a\rrvert ^{2p},\qquad a'
\in\R^q.
\end{eqnarray*}
It is clear that $\psi$ satisfies (\ref{derivpsi}). Moreover, we have by
(\ref{min0un-phi}) and the inf-convolution definition of $u_n$
\begin{eqnarray*}
\psi\bigl(a'\bigr) & \leq& u_n\bigl(t,x,a+a'
- a_n\bigr)-n\llvert a_n-a\rrvert ^{2p} \leq v
\bigl(t,x,a'\bigr),\qquad a ' \in\R^q.
\end{eqnarray*}
Moreover, since $a_n \in\mathring A$ attains the infimum in the
inf-convolution definition of $u_n(t,x,a)$, we have
\begin{eqnarray*}
\psi(a_n) & = & \varphi(a)-n\llvert a_n-a\rrvert
^{2p} = u_n(t,x,a) - n\llvert a_n-a\rrvert
^{2p} = v(t,x,a_n),
\end{eqnarray*}
which shows \eqref{min0upsi}.\vadjust{\goodbreak}

\noindent\textit{Step 3.} \textit{The function $u_n$ does not depend locally
on the variable $a$.}
From Step 2 and Lemma~\ref{lemisolevarasursol}, we obtain that for any
fixed $(t,x) \in[0,T)\times\R^d$, the function $u_n(t,x,\cdot)$
inherits from $v(t,x,\cdot)$ the viscosity supersolution to
%
%
\begin{eqnarray}
\label{viscoun} - \bigl\llvert D_a u_n(t,x,a) \bigr
\rrvert &=& 0,\qquad a \in\mathring A_\delta,
\end{eqnarray}
for any $\delta> 0$, $n \geq N(\llvert x\rrvert,\delta)$.
Let us then show that $u_n(t,x,\cdot)$ is locally constant in the sense
that for all $a \in\mathring{A}_\delta$
%
%
\begin{eqnarray}
\label{consun} u_n(t,x,a) & = & u_n\bigl(t,x,a'
\bigr),\qquad\forall a' \in B(a,\eta),
\end{eqnarray}
for all $\eta>0$ such that $B(a,\eta)\subset\mathring{A}_\delta$. We
first notice from the inf-convolution definition that $u_n(t,x,\cdot)$ is
semiconcave on
$\mathring{A}_\delta$. From Theorem~2.1.7 in \cite{CanSin04}, we deduce
that $u_n(t,x,\cdot)$ is locally Lipschitz continuous on $\mathring
{A}_\delta$. By Rademacher theorem, this implies
that $u_n(t,x,\cdot)$ is differentiable almost everywhere on
$\mathring
{A}_\delta$. Therefore, by Corollary~2.1(ii) in \cite{bar94}, and the
viscosity supersolution property (\ref{viscoun}), we get that
this relation (\ref{viscoun}) holds actually in the classical sense for
almost all $a'\in\mathring{A_\delta}$. In other words,
$u_n(t,x,\cdot)$ is
a locally Lipschitz continuous function with derivatives equal to zero
almost everywhere on $\mathring{A}_\delta$. This means that it is
locally constant (easy exercise in analysis left to the reader).

\noindent\textit{Step 4.} From the convergence (\ref{convun}) of
$u_n$ to $v$,
and the relation (\ref{consun}),
we get by sending $n$ to infinity that for any $\delta>0$ the function
$v$ satisfies: for any $(t,x,a)\in[0,T)\times\R^d\times\mathring
{A}_\delta$
\begin{eqnarray*}
v(t,x,a) & = & v\bigl(t,x,a'\bigr)
\end{eqnarray*}
for all $\eta>0$ such that $B(a,\eta)\subset\mathring{A}_\delta$ and
all $a'\in B(a,\eta)$.
Then by sending $\delta$ to zero we obtain that $v$ does not depend on
the variable $a$ locally on $[0,T)\times\R^d\times\mathring{A}$. Since
$\mathring{A}$ is assumed to be convex,
we obtain that $v$ does not depend on the variable $a$ on $[0,T)\times
\R
^d\times\mathring{A}$.
\end{pf}

\subsection{\texorpdfstring{Viscosity properties of the minimal solution to the BSDE
with $A$-nonposi\-tive jumps}{Viscosity properties of the minimal solution to the BSDE
with $A$-nonpositive jumps}} \label{secpptevisocw}

From Proposition~\ref{THMnondepva}, we can define by misuse of notation
the function $v$ on $[0,T)\times\R^d$ by
%
%
\begin{eqnarray}
\label{vnona} v(t,x) &=& v(t,x,a),\qquad(t,x) \in[0,T)\times\R^d,
\end{eqnarray}
for any $a \in\mathring A$. Moreover, by the growth condition
(\ref{polyvvn}), we have for some $p \geq2$
%
%
\begin{eqnarray}
\label{croislinw} \sup_{(t,x)\in[0,T]\times\R^d} \frac{\llvert v(t,x)\rrvert }{1+\llvert x\rrvert ^p} & < & \infty.
\end{eqnarray}

The aim of this section is to prove that the function $v$ is a
viscosity solution to (\ref{IPDE})--(\ref{termIPDE}).


\begin{pf*}{Proof of the viscosity supersolution property to (\ref{IPDE})}
We first notice from (\ref{vinf1}) and (\ref{vnona}) that $v$ is
l.s.c. and
%
%
\begin{eqnarray}
\label{link-env-inf-w-v} v(t,x) = v_*(t,x) & = & \liminf_{n\rightarrow\infty}{}_*
v_n(t,x,a)
\end{eqnarray}
for all $(t,x,a)\in[0,T]\times\R^d\times\mathring A$. Let $(t,x)$ be
a point in $[0,T)\times\R^d$, and $\varphi\in C^{1,2}([0,T]\times\R
^d)$, such that
\begin{eqnarray*}
(v-\varphi) (t,x) & = & \min_{[0,T]\times\R^d}(v-\varphi).
\end{eqnarray*}
We may assume w.l.o.g. that $\varphi$ satisfies $\sup_{(t,x)\in
[0,T]\times\R^d}\frac{\llvert \varphi(t,x)\rrvert }{1+\llvert x\rrvert ^p} < \infty$.
Fix some $a \in\mathring A$, and define for $\eps> 0$, the
test function
\begin{eqnarray*}
\varphi^\eps\bigl(t',x',a'
\bigr) &=& \varphi\bigl(t',x'\bigr) - \eps \bigl(\bigl
\llvert t'-t\bigr\rrvert ^2 + \bigl\llvert
x'-x\bigr\rrvert ^{2p}+\bigl\llvert a'-a
\bigr\rrvert ^{2p}\bigr),
\end{eqnarray*}
for all $(t',x',a') \in[0,T]\times\R^d\times\R^q$. Since
$\varphi
^\eps(t,x,a) = \varphi(t,x)$, and $\varphi^\eps\leq\varphi$
with equality
iff $(t',x',a') = (t,x,a)$, we then have
%
%
\begin{eqnarray}
\label{strictphir} \bigl(v-\varphi^\eps\bigr) (t,x,a) &=& \mbox{strict}
\min_{[0,T]\times\R
^d\times\R
^q}\bigl(v-\varphi^\eps\bigr).
\end{eqnarray}
From the growth conditions on the continuous functions $v_n$ and
$\varphi$, there exists a bounded sequence $(t_n,x_n,a_n)_n$
(we omit the dependence in $\eps$)
in $[0,T]\times\R^d\times\R^q$ such that
%
%
\begin{eqnarray}
\label{minvnmoinsvarphialpha} \bigl(v_n-\varphi^\eps\bigr)
(t_n,x_n,a_n) & = & \min_{[0,T]\times\R^d\times
\R
^q}
\bigl(v_n-\varphi^\eps\bigr).
\end{eqnarray}
From (\ref{link-env-inf-w-v}) and (\ref{strictphir}), we obtain by
standard arguments that up to a subsequence
\begin{eqnarray*}
\bigl(t_n,x_n,a_n,v_n(t_n,x_n,a_n)
\bigr) & \longrightarrow& \bigl(t,x,a,v(t,x)\bigr), \qquad\mbox{as } n \mbox{ goes
to infinity.}
\end{eqnarray*}
Now from the viscosity supersolution property of $v_n$ at
$(t_n,x_n,a_n)$ with the test function $\varphi^\eps$, we have
\begin{eqnarray*}
&&{-} \frac{\partial\varphi^{\eps}}{\partial t}(t_n,x_n,a_n) - \Lc
^{a_n}\varphi^{\eps
}(t_n,x_n,a_n)
\nonumber
\\
&&{}\qquad-f \bigl(x_n,a_n,v_n(t_n,x_n,a_n),
\sigma\trans (x_n,a_n)D_x\varphi
^{\eps
}(t_n,x_n,a_n),
\Mc^{a_n}\varphi^{\eps}(t_n,x_n,a_n)
\bigr)
\nonumber
\\
&&{}\qquad-\int_A\bigl[\varphi^{\eps}
\bigl(t_n,x_n,a'\bigr)-\varphi^{\eps
}(t_n,x_n,a_n)
\bigr]\lambda _\pi\bigl(da'\bigr) \label{pptesursolvn}
\\
&&{}\qquad- n\int_A\bigl[\varphi^{\eps}
\bigl(t_n,x_n,a'\bigr)-\varphi^{\eps
}(t_n,x_n,a_n)
\bigr]^+\lambda_\pi\bigl(da'\bigr) \geq0 .
\nonumber
\end{eqnarray*}
Sending $n$ to infinity in the above inequality, we get from the
definition of $\varphi^\eps$ and the dominated convergence theorem
%
%
\begin{eqnarray}
\label{ineq-sur-sol-interm} &&{-} \frac{\partial\varphi^{\eps}}{\partial t}(t,x,a) - \Lc ^a
\varphi^{\eps}(t,x,a)
\nonumber
\\
&&\qquad{}-f \bigl(x,a,v(t,x),\sigma\trans(x,a)D_x
\varphi^{\eps
}(t,x,a),\Mc ^a\varphi^{\eps}(t,x,a)
\bigr)
\\
&&\qquad{}+ \eps\int_A\bigl\llvert a'-a
\bigr\rrvert ^{2p} \lambda_\pi\bigl(da'\bigr)
\geq0 .
\nonumber
\end{eqnarray}
Sending $\eps$ to zero, and since $\varphi^\eps(t,x,a) = \varphi
(t,x)$, we get
\begin{eqnarray*}
- \frac{\partial\varphi}{\partial t}(t,x) - \Lc^a\varphi(t,x) -f \bigl(x,a,v(t,x),
\sigma\trans(x,a)D_x\varphi(t,x),\Mc^a\varphi(t,x)
\bigr) & \geq& 0 .
\nonumber
\end{eqnarray*}
Since $a$ is arbitrarily chosen in $\mathring A$, we get from
{(H$A$)} and the continuity of the coefficients $b$, $\sigma$, $\gamma$
and $f$ in the variable $a$
\begin{eqnarray*}
&&{-} \frac{\partial\varphi}{\partial t}(t,x) - \sup_{a\in A} \bigl[
\Lc^a\varphi(t,x) + f \bigl(x,a,v(t,x),\sigma\trans(x,a)D_x
\varphi(t,x),\Mc^a\varphi(t,x) \bigr) \bigr]
\\
&&\qquad\geq0 ,
\end{eqnarray*}
which is the viscosity supersolution property.
\end{pf*}

\begin{pf*}{Proof of the viscosity subsolution property to (\ref{IPDE})}
Since $v$ is the pointwise limit
of the nondecreasing sequence of continuous functions $(v_n)$, and
recalling (\ref{vnona}), we have by \cite{bar94}, page 91:
%
%
\begin{eqnarray}
\label{defv*} v^*(t,x) & = & \limsup_{n\rightarrow\infty}
{}^*v_n(t,x,a)
\end{eqnarray}
for all $(t,x,a)\in[0,T]\times\R^d\times\mathring A$, where
\begin{eqnarray*}
\limsup_{n\rightarrow\infty} {}^*v_n(t,x,a) := \mathop{\mathop {
\limsup_{n\rightarrow\infty}}_{(t',x',a')\rightarrow(t,x,a)}}_{t'<T,
a'\in\mathring A} %
v_n\bigl(t',x',a'\bigr).
\end{eqnarray*}
Fix $(t,x)\in[0,T)\times\R^d$ and $\varphi\in C^{1,2}([0,T]\times
\R
^d)$ such that
%
%
\begin{eqnarray}
\label{v*phi} 
\bigl(v^*-\varphi\bigr) (t,x) & = & \max
_{[0,T]\times\R^d}\bigl(v^*-\varphi\bigr).
\end{eqnarray}
We may assume w.l.o.g. that this maximum is strict and that $\varphi$ satisfies
%
%
\begin{eqnarray}
\label{condcroisvarphi} \sup_{(t,x)\in[0,T]\times\R^d}\frac{\llvert \varphi(t,x)\rrvert }{1+\llvert x\rrvert ^p} & < & \infty.
\end{eqnarray}
%
Fix $a\in\mathring A$ and consider a sequence $(t_n,x_n,a_n)_n$ in
$[0,T)\times\R^d\times\mathring A$ such that
%
%
\begin{eqnarray}
\label{suitean} \bigl(t_n,x_n,a_n,v_n(t_n,x_n,a_n)
\bigr) & \longrightarrow& \bigl(t,x,a,v^*(t,x)\bigr)\qquad\mbox{as } n\rightarrow
\infty.
\end{eqnarray}
Let us define for $n\geq1$ the function $\varphi_n\in
C^{1,2,0}([0,T]\times\R^d\times\R^q)$ by
\begin{eqnarray*}
\varphi_n\bigl(t',x',a'
\bigr) & = & \varphi\bigl(t',x'\bigr)+n \biggl(
{d_{A_{\eta
_n}}(a')\over
\eta_n}\wedge1+\bigl\llvert t'-t_n\bigr
\rrvert ^{2}+\bigl\llvert x'-x_n\bigr\rrvert
^{2p} \biggr),
\end{eqnarray*}
where $A_{\eta_n}$ is defined by (\ref{defAeps}) for $\eps=\eta_n$ and
$(\eta_n)_n$ is a positive sequence converging to 0 s.t.
[such sequence exists by {(H$\lambda_\pi$)}(ii)]
%
%
\begin{eqnarray}
\label{cond2etan} 
n^2\lambda_\pi(A\setminus
A_{\eta_n}) & \longrightarrow& 0\qquad \mbox{as } n\rightarrow\infty.
\end{eqnarray}
From the growth conditions (\ref{croislinw}) and (\ref{condcroisvarphi})
on $v$ and $\varphi$, we can find a sequence $(\bar t_n,\bar x_n,\bar
a_n)$ in $[0,T]\times\R^d\times A$ such that
%
%
\begin{eqnarray}
\label{minvnmoinsvarphialpha2} (v_n-\varphi_n) (\bar t_n,
\bar x_n,\bar a_n) & = & \max_{[0,T]\times
\R
^d\times A}(v_n-
\varphi_n),\qquad n\geq1 .
\end{eqnarray}
%
Using (\ref{defv*}) and (\ref{v*phi}), we obtain by standard arguments
that up to a subsequence
%
%
\begin{eqnarray}
\label{conv suiteoptimisante2} 
\qquad n \biggl({1\over\eta_n}d_{A_{\eta_n}}(
\bar a_n)+\llvert \bar t_n -t_n\rrvert
^{p}+\llvert \bar x_n-x_n\rrvert
^{2p} \biggr) & \longrightarrow& 0\qquad\mbox{as } n \rightarrow\infty,
\end{eqnarray}
and
\begin{eqnarray*}
v_n(\bar t_n,\bar x_n,\bar a_n
) & \longrightarrow& v^*(t,x)\qquad \mbox{as } n \rightarrow\infty.
\end{eqnarray*}
We deduce from (\ref{conv suiteoptimisante2}) and (\ref{suitean}) that,
up to a subsequence
%
%
\begin{eqnarray}
\label{conv suiteoptimisante3} (\bar t_n,\bar x_n,\bar a_n)
& \longrightarrow& (t,x, \bar a),\qquad \mbox{as } n \rightarrow\infty,
\end{eqnarray}
for some $\bar a \in A$. Moreover, for $n$ large enough, we have
$\bar a_n\in\mathring A$. Indeed, suppose that, up to a subsequence,
$\bar a_n\in\partial A $ for $n\geq1$.
Then we have ${1\over\eta_n}d_{A_{\eta_n}}(\bar a_n)\geq1$, which
contradicts (\ref{conv suiteoptimisante2}).
Now, from the viscosity subsolution property of $v_n$ at $(\bar
t_n,\bar x_n,\bar a_n)$ with the test function $\varphi_n$ satisfying
(\ref{minvnmoinsvarphialpha2}), Lemma~\ref{lemviscovnA}, and since
$\bar
a_n\in\mathring A$, we have
%
\begin{eqnarray}
\label{inegphin} &&{-} \frac{\partial\varphi_n}{\partial t}(\bar t_n,\bar
x_n,\bar a_n) - \Lc^{\bar a_n}\varphi _n(
\bar t_n,\bar x_n,\bar a_n)
\nonumber
\\
&&{}\qquad-f\bigl(\bar x_n,\bar a_n,v_n(\bar
t_n,\bar x_n,\bar a_n),\sigma \trans (\bar
x_n,\bar a_n)D_x\varphi(\bar
t_n,\bar x_n),\Mc^{\bar a_n}\varphi (\bar
t_n,\bar x_n,\bar a_n)\bigr)
\\
&&{}\qquad-(n+1)n\int_A \biggl({d_{A_{\eta_n}}(a')\over\eta_n}
\wedge 1 \biggr)\lambda _\pi\bigl(da'\bigr) \leq 0 ,
\nonumber
\end{eqnarray}
for all $n\geq1$. From (\ref{cond2etan}), we get
%
%
\begin{eqnarray}
\label{cvtermintsssol} (n+1)n\int_A \biggl({d_{A_{\eta_n}}(a')\over\eta_n}
\wedge1 \biggr)\lambda_\pi \bigl(da'\bigr)&
\longrightarrow& 0\qquad\mbox{as } n\rightarrow\infty.
\end{eqnarray}
Sending $n$ to infinity into (\ref{inegphin}), and using (\ref{defv*}),
(\ref{conv suiteoptimisante3}) and (\ref{cvtermintsssol}), we get
\begin{eqnarray*}
- \frac{\partial\varphi}{\partial t}(t,x) - \Lc^{\bar a}\varphi (t,x) 
-f
\bigl(x,\bar a,v^*(t,x),\sigma\trans(x,\bar a)D_x\varphi(t,x),\Mc
^{\bar
a}\varphi(t,x)\bigr) & \leq& 0 .
\end{eqnarray*}
Since $\bar a \in A$, this gives
\begin{eqnarray*}
&&{-} \frac{\partial\varphi}{\partial t}(t,x) - \sup_{a\in A} \bigl[
\Lc^a\varphi(t,x) + f \bigl(x,a,v^*(t,x),\sigma\trans(x,a)D_x
\varphi(t,x),\Mc^a\varphi(t,x) \bigr) \bigr]
\\
&&\qquad\leq0 ,
\end{eqnarray*}
which is the viscosity subsolution property.
\end{pf*}

\begin{pf*}{Proof of the viscosity supersolution property to (\ref{termIPDE})}
Let $(x,a)\in\mathbb{R}^d\times\mathring A$. From (\ref
{link-env-inf-w-v}), we can find a sequence $(t_n,x_n,a_n)_n$ valued in
$[0,T)\times\R^d\times\R^q$ such that
\begin{eqnarray*}
\bigl(t_n,x_n,a_n,v_n(t_n,x_n,a_n)
\bigr) & \longrightarrow& \bigl(T,x,a,v_*(T,x)\bigr)\qquad\mbox{as } n \rightarrow
\infty.
\end{eqnarray*}
The sequence of continuous functions $(v_{n})_{n}$ being nondecreasing
and $v_{n}(T,\cdot)=g$ we have
\begin{eqnarray*}
v_*(T,x) & \geq& \lim_{n\rightarrow\infty}v_1(t_n,x_n,a_n)
= g(x,a).
\end{eqnarray*}
Since $a$ is arbitrarily chosen in $\mathring A$, we deduce that
$v_*(T,x) \geq\sup_{a\in\mathring A}g(x,a) = \sup_{a\in A}
g(x,a)$ by {(H$A$)} and continuity of $g$ in $a$.
\end{pf*}

\begin{pf*}{Proof of the viscosity subsolution property to (\ref{termIPDE})}
Let $x\in\mathbb{R}^d$. Then we can find by (\ref{defv*})
a sequence $(t_n,x_n,a_n)_n$ in $[0,T)\times\R^d\times\mathring A$
such that
%
%
\begin{eqnarray}
\label{recupw^*} \bigl(t_n,x_n,v_n(t_n,x_n,a_n)
\bigr) & \rightarrow& \bigl(T,x,v^*(T,x)\bigr),\qquad\mbox {as } n\rightarrow
\infty.
\end{eqnarray}
Define the function $h\dvtx[0,T]\times\R^d\rightarrow\R$ by
\begin{eqnarray*}
h(t,x) & = & \sqrt{T-t}+\sup_{a\in A}g(x,a)
\end{eqnarray*}
for all $(t,x)\in[0,T)\times\R^d$. From {(HFC)}, {(HBC1)}
and {(HBC2$'$)}, we see that $h$ is a continuous viscosity
supersolution to (\ref{IPDEn})--(\ref{termIPDEn}), on $[T-\eta,T]\times
\bar B(x,\eta)$ for $\eta$ small enough. We can then apply Theorem~3.5
in \cite{barbucpar97} which gives that
\begin{eqnarray*}
\label{majv_nh} v_n & \leq& h\qquad\mbox{on } [T-\eta,T]\times\bar
B(x,\eta )\times A
\end{eqnarray*}
for all $n\geq0$. By applying the above inequality at
$(t_n,x_n,a_n)$, and sending $n$ to infinity, together with (\ref
{recupw^*}), we get the required result.
\end{pf*}

\section{Conclusion}\label{sec4}

We introduced a class of BSDEs with partially nonpositive jumps and
showed how the minimal solution is related to a fully nonlinear IPDE of
HJB type, when considering a Markovian framework with forward regime
switching jump-diffusion process.
Such BSDE representation can be extended to cover also the case of
the Hamilton--Jacobi--Bellman--Isaacs equation arising in
controller/stopper differential games; see \cite{chocospha13}.
It is also extended to the non-Markovian context in \cite{fuhpha13}.
From a numerical application viewpoint, our BSDE representation leads
to original probabilistic approximation scheme for the resolution in
high dimension of fully nonlinear HJB equations, as recently
investigated in \cite{KLP13a} and \cite{KLP13b}. We believe that this
opens new perspectives for dealing with more general nonlinear
equations and control problems, like for instance mean field games or
control of McKean--Vlasov equations.


\section*{Acknowledgments}
The authors would like to thank Pierre Cardaliaguet for useful discussions.



\printaddresses

\end{document}